\documentclass[12pt,reqno]{amsart}
\usepackage{hyperref}
\hypersetup{
    colorlinks=true,       % false: boxed links; true: colored links
    linkcolor=blue,          % color of internal links
    citecolor=blue,        % color of links to bibliography
    filecolor=blue,      % color of file links
    urlcolor=cyan  
}
\usepackage{latexsym,amssymb,amsthm,amsmath,mathtools}

\newcommand{\R}{{\ensuremath{\mathbb{R}}}}
\newcommand{\N}{{\ensuremath{\mathbb{N}}}}

\renewcommand{\P}{\ensuremath{\mathbb{P}}}
\renewcommand{\dj}{d\kern-0.4em\char"16\kern-0.1em}
\newcommand{\E}{\ensuremath{\mathbb{E}}}

\newtheorem{Thm}{Theorem}[section]

\newtheorem{Cor}[Thm]{Corollary}
\newtheorem{Lem}[Thm]{Lemma}
\newtheorem{Prop}[Thm]{Proposition}

\theoremstyle{remark}
\newtheorem{Rem}[Thm]{Remark}

\theoremstyle{definition}
\newtheorem{Ex}[Thm]{Example}

\theoremstyle{definition}

\theoremstyle{definition}

\begin{document}
\numberwithin{equation}{section}
\bibliographystyle{amsalpha}

\title[Asymptotical properties of distributions]{Asymptotical properties of distributions of isotropic L\' evy processes}
\begin{abstract}
In this paper, we establish the precise asymptotic behaviors of the tail probability and the transition density of a large class of isotropic L\'evy processes when the scaling order  is between 0 and 2 {\it including 2}. We also obtain the precise asymptotic behaviors of the tail probability of subordinators when the scaling order  is between 0 and 1 {\it including 1}.

 The asymptotic
 expressions are given in terms of the radial part of characteristic exponent  $\psi$ and its derivative. In particular, when  $\psi(\lambda)-\frac{\lambda}{2}\psi'(\lambda)$ varies regularly, as  $\frac{t\psi(r^{-1})^2}{\psi(r^{-1})-(2r)^{-1}\psi'(r^{-1})} \to 0$ the tail probability  $\P(|X_t|\geq r)$ is asymptotically equal to  a constant times $ t( \psi(r^{-1})-(2r)^{-1}\psi'(r^{-1})).$
\end{abstract}

\author{Panki Kim}
\address{Department of Mathematical Sciences,
Seoul National University,
Building 27, 1 Gwanak-ro, Gwanak-gu
Seoul 08826, Republic of Korea}
\curraddr{}
\thanks{This work was supported by the National Research Foundation of Korea(NRF) grant funded by the Korea government(MSIP) 
(No. 2016R1E1A1A01941893).}
\email{pkim@snu.ac.kr}
\author[A.\ Mimica]{Ante Mimica\,$(\dagger)$}
 %\url{https://web.math.pmf.unizg.hr/~amimica/ }
% \address{Department of Mathematics, University of Zagreb, Bijeni\v cka cesta 30, 10 000 Zagreb, Croatia}
% \curraddr{ \url{https://web.math.pmf.unizg.hr/~amimica/ }}
%% \thanks{Supported in part by Croatian Science Foundation under the project 3526.}
%% \email{amimica@math.hr}
\address[Mimica]{Ante Mimica, *20-Jan-1981\,-- %(Split, Croatia),
	$\dagger$\,9-Jun-2016 %(Zagreb, Croatia)
	\\
\url{https://web.math.pmf.unizg.hr/~amimica/}}



\maketitle

\allowdisplaybreaks[3]
\section{Introduction and main results}

This paper is a continuation of the journey of the second-named author on investigating the behavior of the transition density of  L\'evy processes whose (weak) scaling order  is between 0 and 2 {\it including 2}.
In \cite{Mi6}, the second-named author considered the large class of
purely  discontinuous subordinate Brownian motions when the  weak scaling order  is between 0 and 2  including 2, and obtained sharp heat kernel estimates of such processes. 
In this paper, we assume that 
the function $\psi(\lambda)-\frac{\lambda}{2}\psi'(\lambda)$ 
is regularly varying where $\xi \to \psi(|\xi|)$ is  the characteristic exponent of an isotropic L\'evy process, and 
discuss 
the asymptotic behavior of the tail probability. 
As a corollary,  when the process is unimodal,   
the asymptotic behavior of the tail probability implies the precise off-diagonal asymptotic expression of the transition density.
We remark here that 
the function $\psi(\lambda)-\frac{\lambda}{2}\psi'(\lambda)$  is a natural counterpart of the function 
$H(\lambda)=\phi(\lambda)-\lambda\phi'(\lambda)$ 
 for subordinator with Laplace exponent $\phi$, which 
has been already known in the literature (see \cite{JP}).

L\'evy  processes and their associated (non-local) operators
have been of current research interest both in probability theory and in PDE. 
 The transition density  
 of a L\'evy  process $X$ is also the fundamental solution  of corresponding non-local operator.
 Except a few special case,  the explicit expression of transition density of L\'evy  process is
  typically impossible to get. 
Thus obtaining the exact  off-diagonal asymptotic expression of transition density is 
an important problem both in probability theory and in analysis. 

Even though transition densities of L\'evy processes are  determined by their characteristic exponents, in general, 
neither the transition density nor the tail probability can be calculated explicitly by their characteristic exponents. 
But for large class of isotropic L\'evy process, their  asymptotic behaviors can be formulated in terms of their characteristic exponents. For example, for the transition density $p(t,x)$  of isotropic $\alpha$-stable process, whose characteristic exponent is $|x|^\alpha$, it is known that
\begin{align}
\lim_{  t |x|^{-\alpha}\to 0 }  \frac{p(t,x)}{ t|x|^{-d-\alpha}}  =\alpha2^{\alpha-1}\pi^{-d/2-1}\sin\left(\frac{\alpha\pi}{2}\right)\Gamma\left(\frac{\alpha}{2}\right)\Gamma\left(\frac{d+\alpha}{2}\right),\label{e:BG}
\end{align}
where $\Gamma(t)$ is the gamma function. (See \cite{BG}.)
Since, in general,  the L\'evy kernel of  isotropic L\'evy process $X$ is not equal to a constant times $|x|^{-d}\psi(|x|^{-1})$
where $\xi \to \psi(|\xi|)$ is  the characteristic exponent of $X$, we do not expect that \eqref{e:BG} holds. 
But as it is shown in \cite[Theorems 4 and 5]{CGT}, when  
 the characteristic exponent $\psi$ of isotropic unimodal L\'evy process
 is regularly varying at $\infty$ (at zero, respectively)  with index $\alpha \in (0,2)$, 
the transition density is asymptotically equal to a constant times 
$t|x|^{-d}\psi(|x|^{-1})$ as $t \psi(|x|^{-1})\to 0$ and $|x| \to \infty$ ($|x| \to 0$, respectively).

Natural open questions are, when  the regularly varying index $\alpha$ is not strictly between $0$ and $2$, 

\noindent
(1)  instead of $t \psi(|x|^{-1})\to 0$, what condition in the limit we need to put?

\noindent
(2) to what functions the tail probability and the transition density, respectively, are asymptotically equal? 

In this paper, we answer the above questions when  $\psi(\lambda)-\frac{\lambda}{2}\psi'(\lambda)$ varies regularly.
Up to our knowledge,  the asymptotic expressions in 
Theorems \ref{thm:subord_tail_estimate}, \ref{t:tail}(ii) and \ref{cor:conseq2}(ii) are new forms comparing to existing results (for example, \cite{Be, BG, CGT, FJH, Is, Kn, Le}) and
 it also covers \cite[Theorems 4 and 5]{CGT}, see Theorems \ref{t:tail}(i) and \ref{cor:conseq2}(i).

In this paper, 
$f(t) \sim g(t)$, $t \to a$ means $
\lim_{t \to a} f(t)/g(t) = 1$. 
We also  use the notation
$f(t, r) \sim g(t,r)$, $r \to a$ and  $h(t,r) \to 0$,
which means 
 $$
\lim_{r \to a \atop h(t,r) \to 0 } \frac{f(t,r)}{g(t,r)} = 1.$$ 

We say a function $\ell : (0, \infty) \to (0, \infty)$
varies regularly at $\infty$ (at zero, respectively) with index $\alpha$ if 
 $\ell(\lambda t)/\ell(t) \to \lambda^\alpha$ as $t \to \infty$ ($t \to 0$, respectively) for every $ \lambda>0$. 

$T=(T_t)_{t\geq 0}$ is called a  subordinator if it is a non-decreasing L\' evy process. The Laplace transform of the law of a  subordinator  $T_t$ is given by
\begin{equation}\label{eq:lap_tf}
    \E e^{-\lambda T_t}=e^{-t\phi(\lambda)},\quad \lambda>0\,,
\end{equation}
and  $\phi:(0,\infty)\rightarrow (0,\infty)$ is called the Laplace exponent of the subordinator $T$.

The first main result of this paper is the following asymptotic expression of the tail probability of  a subordinator. 
 \begin{Thm}\label{thm:subord_tail_estimate}
 Suppose that $\phi$ is the Laplace exponent of a subordinator $T$.
If $H(\lambda)=\phi(\lambda)-\lambda\phi'(\lambda)$ varies regularly at $0$ (at $\infty$, respectively) with index 
$\gamma\in [0,2)$, then 
\begin{align}
\label{e:t.1.1}
    {\P(T_t\geq r)}{} \sim\frac{1}{\Gamma(2-\gamma)}  t H(r^{-1})
, 
   \quad  r\to\infty \,\, (r \to 0, \text{ resp.})\, \text{ and }\, t\frac{\phi(r^{-1})^2}{H(r^{-1})} \to 0
\end{align}
 \end{Thm}
 
 It turns out  that, in fact, the tail of the L\' evy measure of $T$ describes decay of $\P(T_t\geq r)$.
(See Proposition \ref{prop:levy-tail} below.)

Note that, 
if $\phi$ varies regularly at $0$ (at $\infty$, respectively) with index 
$\gamma \in (0,1)$, 
since $\phi'$ is non-increasing, we can apply monotone density theorem (see \cite[Theorem 1.7.2]{BGT}) to get 
$\phi'(\lambda)\sim \gamma \lambda^{-1}\phi(\lambda)$ 
as $\lambda \to 0$. Hence,
\begin{equation}\label{eq:H_phi}
\lim_{\lambda\to 0\atop
(\lambda\to \infty, \text{resp.})}\frac{H(\lambda)}{\phi(\lambda)}=1-\lim_{\lambda\to 0\atop
(\lambda\to \infty, \text{resp.})}\frac{\lambda \phi'(\lambda)}{\phi(\lambda)}=1-\gamma\,.
\end{equation}
\eqref{eq:H_phi} and  Proposition \ref{prop:Hregvarn} below imply that 
 if $\gamma\in (0,1)$  then $\phi$ varies regularly at $0$ with index $\gamma$ if and only if  $H$  varies regularly at $0$  with index $\gamma$. Furthermore, from \eqref{eq:H_phi} we see that for $\gamma=1$ functions $H$ and $\phi$ are not comparable, but we will see that regular variation of $H$ implies regular variation of $\phi$ (see Proposition \ref{prop:Hregvarn} below).

A consequence of  Theorem \ref{thm:subord_tail_estimate} and \eqref{eq:H_phi} 
is that, if $\gamma\in [0,1)$, we may replace condition $t\frac{\phi(r^{-1})^2}{H(r^{-1})}\to 0$ in \eqref{e:t.1.1} by $t\phi(r^{-1})\to 0$.
\begin{Cor}\label{cor:asympt}
Assume that $\phi$ varies regularly at $0$ (at $\infty$, respectively) with index $\gamma\in [0,1)$. Then
 \[
    {\P(T_t\geq r)}{} \sim\frac{1}{\Gamma(1-\gamma)}   t\phi(r^{-1})
, 
   \quad  r\to\infty \,\, (r \to 0, \text{ resp.}) \quad \text{ and } \quad  t \phi(r^{-1}) \to 0.
 \]
\end{Cor}

The next main result is 
 the precise asymptotic expression of the tail probability of isotropic L\'evy processes.
We recall that $\xi \to \psi(|\xi|)$ is the characteristic exponent  of isotropic L\'evy process $X$ if 
\[
    \E[e^{i\xi\cdot X_t}]=e^{-t\psi(|\xi|)},\quad \xi\in \R^d\,,
\]
where 
\begin{align}
\label{e:psirep}
    \psi(|\xi|)=a|\xi|^2+\int_{\R^d}(1-e^{i\xi\cdot x})\,\nu(dx)
\end{align}
 with $a\geq 0$ and a L\' evy measure $\nu$\,.

\begin{Thm}\label{t:tail}
Suppose $X=(X_t)_{t\geq 0}$ is an isotropic L\' evy process in $\R^d$, $d\geq 1$ with the characteristic exponent $\xi \to \psi(|\xi|)$.
Let $g(\lambda)=\psi(\lambda)-\frac{\lambda}{2}\psi'(\lambda)$.
%\noindent
%(i) If $\psi$ varies regularly at $0$ (at $\infty$, respectively) with index $\alpha\in [0,2)$, then 
%\begin{align*}
%\P(|X_t|\geq r)  &\sim \frac{(2-\alpha)2^{\alpha-1}\Gamma(\frac{d+\alpha}{2})}{\Gamma(2-\frac{\alpha}{2})\Gamma(\frac{d}{2})} t\psi(r^{-1}), 
%\\
% &\quad  r\to\infty \,\, (r \to 0, \text{ resp.})
% \quad \text{ and }  \quad
%  t \psi(r^{-1})\to 0\,.
%\end{align*}

\noindent
(i) Suppose that  $g$ varies regularly at $0$ with index $\alpha\in 
[0, 4)$.  Then 
\begin{align*}
\P(|X_t|\geq r) &\sim
 \frac{2^{\alpha-1}\Gamma(\frac{d+\alpha}{2})}{\Gamma(2-\frac{\alpha}{2})\Gamma(\frac{d}{2})} tg(r^{-1}),
\quad  r\to\infty\\
 &\ 
 \quad \text{ and }  \quad
 \frac{t\psi(r^{-1})^2}{\psi(r^{-1})-(2r)^{-1}\psi'(r^{-1})}=  \frac{t\psi(r^{-1})}{1-\frac{r^{-1}\psi'(r^{-1})}{2\psi(r^{-1})}}\to 0\,.
\end{align*}

\noindent
(ii) Suppose that  $g$ varies regularly at $\infty$ with index $\alpha \ge 0$. Then, in fact,  $g$ varies regularly at $\infty$ with index $\alpha \in [0, 2]$.
We further assume that $\psi$ does not have diffusion part ($a=0$ in \eqref{e:psirep}) if $g$ varies regularly at $\infty$ with index $\alpha\in 
[0, 2)$.
 Then 
\begin{align*}
\P(|X_t|\geq r) &\sim
 \frac{2^{\alpha-1}\Gamma(\frac{d+\alpha}{2})}{\Gamma(2-\frac{\alpha}{2})\Gamma(\frac{d}{2})} tg(r^{-1}),
\quad  r\to0 \text{ and } 
 \frac{t\psi(r^{-1})^2}{\psi(r^{-1})-(2r)^{-1}\psi'(r^{-1})}\to 0\,.
\end{align*}
\end{Thm}

Finally we give 
 the precise off-diagonal asymptotic expression of  the transition density $p(t,x)$ when the isotropic L\'evy process is  unimodal, that is, 
 $r \to q(t,r)$ is decreasing where $p(t,x)=q(t,|x|)$. 

\begin{Thm}\label{cor:conseq2}
Suppose $X=(X_t)_{t\geq 0}$ is  a unimodal L\' evy process in $\R^d$, $d\geq 1$ with the characteristic exponent $\xi \to \psi(|\xi|)$, and that 
$p(t,x)$ and $J(x)$ are the transition density and the L\'evy  kernel of $X$, respectively. 
%\noindent
%(i) If $\psi$ varies regularly at $0$ (at $\infty$, respectively) with index $\alpha\in (0,2)$, then 
%\begin{align}\label{e:P1}
%p(t,x)  &\sim \alpha 2^{\alpha-1}\pi^{-d/2}\frac{\Gamma(\frac{d+\alpha}{2})}{\Gamma(1-\frac{\alpha}{2})}  t|x|^{-d}\psi(|x|^{-1})
%\nonumber\\&=\alpha2^{\alpha-1}\pi^{-d/2-1}\sin\left(\frac{\alpha\pi}{2}\right)\Gamma\left(\frac{\alpha}{2}\right)\Gamma\left(\frac{d+\alpha}{2}\right)  t|x|^{-d}\psi(|x|^{-1}),
%\nonumber\\\
% &\quad  |x|\to\infty \,\, (|x| \to 0, \text{ resp.})
% \quad \text{ and }  \quad
%  t \psi(|x|^{-1})\to 0\,,
%\end{align}
%and \begin{align}\label{e:J1}
%J(x)  \sim
%\alpha 2^{\alpha-1}\pi^{-d/2}\frac{\Gamma\left(\frac{d+\alpha}{2}\right)}{\Gamma\left(1-\frac{\alpha}{2}\right)}|x|^{-d}\psi(|x|^{-1}),
%\quad  |x|\to\infty \,\, (|x| \to 0, \text{ resp.}). \end{align}
Suppose that  $g(\lambda)=\psi(\lambda)-\frac{\lambda}{2}\psi'(\lambda)$ varies regularly at $0$ with index $\alpha\in 
(0, 4)$  (at $\infty$ with index $\alpha\in 
(0, 2]$, respectively). We further assume that $\psi$ does not have diffusion part  ($a=0$ in \eqref{e:psirep}) if $g$ varies regularly at $\infty$ with index $\alpha\in 
(0,  2)$.
 Then 
\begin{align}\label{e:P2}
p(t,x)  &\sim
\alpha 2^{\alpha-1}\pi^{-d/2}\frac{\Gamma\left(\frac{d+\alpha}{2}\right)}{\Gamma\left(2-\frac{\alpha}{2}\right)}t|x|^{-d}g(|x|^{-1}),
\nonumber\\\
 &\quad  |x|\to\infty \,\, (|x| \to 0, \text{ resp.})
 \quad \text{ and }  \quad
  \frac{t\psi(|x|^{-1})}{1-\frac{|x|^{-1}\psi'(|x|^{-1})}{2\psi(|x|^{-1})}}\to 0\,,
\end{align}
and \begin{align}\label{e:J2}
J(x)  &\sim
\alpha 2^{\alpha-1}\pi^{-d/2}\frac{\Gamma\left(\frac{d+\alpha}{2}\right)}{\Gamma\left(2-\frac{\alpha}{2}\right)}|x|^{-d}g(|x|^{-1}),
\quad  |x|\to\infty \,\, (|x| \to 0, \text{ resp.}). \end{align}
\end{Thm}

We conclude this introduction by setting up some notation and conventions.
We use ``$:=$" to denote a definition, which
is  read as ``is defined to be"; we denote $a \wedge b := \min \{ a, b\}$,
$a \vee b := \max \{ a, b\}$.

\section{Tail asymptotic behaviors of subordinators}

In this section we always assume that $T=(T_t)_{t\geq 0}$ is a subordinator and that  $\phi:(0,\infty)\rightarrow (0,\infty)$ is the Laplace exponent of $T$.
$\phi$ has the following form (see \cite[Section 3]{Be})
\begin{equation}\label{eq:lap_exp}
    \phi(\lambda)=b\lambda +\int_{(0,\infty)}(1-e^{-\lambda t})\mu(dt)\,.
\end{equation}
Here $b\geq 0$ is called the drift  and $\mu$ is the L\' evy measure of the subordinator $T$, i.e. a measure on $(0,\infty)$ satisfying $\int_{(0,\infty)}(1\wedge t)\mu(dt)<\infty$\,.
Laplace exponent $\phi$ belongs to the class of Bernstein functions, i.e., a non-negative $C^\infty(0,\infty)$ function such that  $(-1)^{n}\phi^{(n)}\leq 0$ for all $n\in \N$.

Let us define the function $H:(0,\infty)\rightarrow [0,\infty)$ by 
\begin{align}
\label{e:dH1}
H(\lambda)=\phi(\lambda)-\lambda\phi'(\lambda), \quad \lambda>0.
\end{align} 

Note that, by the concavity of $\phi$,  $H(\lambda) \ge 0$. 
Moreover, $H$ is non-decreasing since $H'(\lambda)=-\lambda \phi''(\lambda) \ge 0$.
We also note that 
\begin{align}
H(\lambda)&=-\lambda^2 (\frac{\phi(\lambda)}{\lambda})' \label{e:dH3}.
\end{align} 
We remark  here that $H$ loses the information on the drift of $\phi$.

In this section we obtain asymptotical properties of the function $\P(T_t\geq r)$ as $r\rightarrow \infty$ ($r\rightarrow 0$, respectively) and $t\frac{\phi(r^{-1})^2}{H(r^{-1})}\rightarrow 0$\,. Note, in particular, that in this case $t\phi(r^{-1})\rightarrow 0$,  since $H(\lambda)\leq \phi(\lambda)$ for all $\lambda>0$\,.

Using \eqref{eq:lap_tf} and Fubini theorem, we obtain 
\begin{equation}\label{eq:tail_lt}
\int_0^\infty e^{-\lambda r}\P(T_t\geq r)\,dr=\E \int_0^{T_t} e^{-\lambda r}\,dr=\E\left[\frac{1-e^{-\lambda T_t}}{\lambda}\right]=\frac{1-e^{-t\phi(\lambda)}}{\lambda}\,.
\end{equation}
 For a measure $\nu$ on $[0,\infty)$ its Laplace transform $\mathcal{L}\nu$ is defined by 
 \[
     (\mathcal{L}\nu)(\lambda)=\int_{[0,\infty)} e^{-\lambda y}\nu(dy)\,, \quad \lambda>0\,.
 \]
 In particular, the Laplace transform of the measure $\nu(dy)=\P(T_t\geq y)\,dy$ is, by \eqref{eq:tail_lt}, given by $(\mathcal{L}\nu)(\lambda)=\frac{1-e^{-t\phi(\lambda)}}{\lambda}$ for $\lambda>0$\,.
 
  We will use the following result which is a part of the continuity theorem for Laplace transforms (see \cite[Theorem XIII.1.2]{Fel} for a proof).
 \begin{Thm}\label{thm:thm_cont}
 Let $(\nu_n)_{n\in \N}$ be a sequence of measures on $[0,\infty)$ and $\lambda_0\geq 0$ such that $\Phi(\lambda):=\lim\limits_{n\to\infty} (\mathcal{L}\nu_n) (\lambda)$ exists for all $\lambda>\lambda_0$\,. Then $\Phi$ is the Laplace transform of a measure $\nu$ on $[0,\infty)$ and the following convergence holds
 \[
     \lim_{n\to\infty} \nu_n([0,x])=\nu([0,x])
 \]
 for all $x\geq 0$ such that $y\mapsto \nu([0,y])$ is continuous at $y$\,.
 \end{Thm}
  
 Now we prove the main result of this section, namely the tail estimate of the subordinator. 

\proof[Proof of Theorem \ref{thm:subord_tail_estimate}]
 Taking derivative in \eqref{eq:tail_lt} we obtain
 \begin{equation}\label{eq:tail_lt_der}
 \int_0^\infty e^{-\lambda r}r\P(T_t\geq r)\,dr=\frac{1-e^{-t\phi(\lambda)}-t\lambda\phi'(\lambda)e^{-t\phi(\lambda)}}{\lambda^2}\,, \quad \lambda>0\,.
 \end{equation}
Using 
\begin{align*}
 & 1-e^{-t\phi(\lambda x^{-1})}-t\lambda x^{-1}\phi'(\lambda x^{-1})e^{-t\phi(\lambda x^{-1})}\\
 =&1-e^{-t\phi(\lambda x^{-1})}-t\phi(\lambda x^{-1})e^{-t\phi(\lambda x^{-1})}+tH(\lambda x^{-1})e^{-t\phi(\lambda x^{-1})}
\end{align*}

and 
\[
    1-e^{-t\phi(\lambda x^{-1})}
    -t\phi(\lambda x^{-1})e^{-t\phi(\lambda x^{-1})}\leq \frac{1}{2}(t\phi(\lambda x^{-1}))^2
\]    
in \eqref{eq:tail_lt_der} it follows that
  \begin{align}\label{eq:tail_lt_der0}
&\frac{ \int_0^\infty e^{-\lambda x^{-1} r}r\P(T_t\geq r)\,dr}{tx^2H(\lambda x^{-1})}=\frac{1-e^{-t\phi(\lambda x^{-1})}-t\lambda x^{-1}\phi'(\lambda x^{-1})e^{-t\phi(\lambda x^{-1})}}{\lambda^2  {tH(\lambda x^{-1})}}\nonumber\\
=& \frac{1-e^{-t\phi(\lambda x^{-1})}-t\phi(\lambda x^{-1})e^{-t\phi(\lambda x^{-1})}}{\lambda^2  {tH(\lambda x^{-1})}}+{\lambda^{-2 } }{e^{-t\phi(\lambda x^{-1})}}\nonumber\\
\le& \frac{t \phi(\lambda x^{-1})^2}{2\lambda^2  {H(\lambda x^{-1})}}+{\lambda^{-2 } }{e^{-t\phi(\lambda x^{-1})}}\,,
 \quad \lambda, x>0\,.
 \end{align}
Note that, by using $\phi(\lambda y)\leq (1\vee y)\phi(\lambda)$ and \cite[Lemma 2.1]{Mi6},  
    \[
        0\leq t\frac{\phi(\lambda x^{-1})^2}{H(\lambda x^{-1})}=t\frac{\phi(x^{-1})^2}{H(x^{-1})}\left(\frac{\phi(\lambda x^{-1})}{\phi(x^{-1})}\right)^2\frac{H( x^{-1})}{H(\lambda x^{-1})}\leq t\frac{\phi(x^{-1})^2}{H(x^{-1})}(1\vee \lambda^2)(1\vee \lambda^{-2})\,.
    \]
Thus 
  $t\frac{\phi(x^{-1})^2}{H(x^{-1})}\to 0$ 
    implies that $t\frac{\phi(\lambda x^{-1})^2}{H(\lambda x^{-1})}\to 0$ for any $\lambda>0$ (and so $t\phi(\lambda x^{-1})\to 0$ for any $\lambda>0$). Therefore from \eqref{eq:tail_lt_der0} we see that 
    \[
    \lim_{ t\frac{\phi(x^{-1})^2}{H(x^{-1})}\to 0}\frac{\int_0^\infty e^{-\lambda x^{-1} r}r \P(T_t\geq r)\,dr}{tx^2H(\lambda x^{-1})}=\lambda^{-2}\quad \text{ for all }\,\,\, \lambda>0\,.
    \]
We now assume that $H(\lambda)=\phi(\lambda)-\lambda\phi'(\lambda)$ varies regularly at $0$ (at $\infty$, respectively) with index 
$\gamma\in [0,2)$. Then, 
  changing variables in the integral of the last display and using regular variation of $H$, we conclude that for any $\lambda>0$ the following holds
\begin{align*}
\lambda^{-2}&=\lim_{x\to\infty 
    (x \to 0, \text{ resp.})  
     \atop t\frac{\phi(x^{-1})^2}{
H(x^{-1})}\to 0}\frac{\int_0^\infty e^{-\lambda r}xr \P(T_t\geq xr)x\,dr}{tx^2H(x^{-1})}\cdot \frac{H(x^{-1})}{H(\lambda x^{-1})}\\&=\lim_{x\to\infty 
    (x \to 0, \text{ resp.})  
     \atop t\frac{\phi(x^{-1})^2}{
H(x^{-1})}\to 0}\frac{\int_0^\infty e^{-\lambda r}r \P(T_t\geq xr)\,dr}{tH(x^{-1})}\lambda^{-\gamma}\,.
\end{align*}
We rewrite the last display  in the following way;
\begin{align}\label{eq:lt-tail-distr}
\lim_{x\to\infty 
    (x \to 0, \text{ resp.})  
     \atop t\frac{\phi(x^{-1})^2}{
H(x^{-1})}\to 0}
\frac{\int_0^\infty e^{-\lambda r}r \P(T_t\geq xr)\,dr}{tH(x^{-1})}=\lambda^{\gamma-2}=\int_0^\infty e^{-\lambda r}\frac{r^{1-\gamma}}{\Gamma(2-\gamma)}\,dr\,.
\end{align}
Hence, by using Theorem \ref{thm:thm_cont} we obtain that for all $R>0$ 
\[
    \lim_{x\to\infty 
    (x \to 0, \text{ resp.})  
     \atop t\frac{\phi(x^{-1})^2}{
H(x^{-1})}\to 0} \frac{\int_0^R r\P(T_t\geq xr)\,dr}{tH(x^{-1})}=\int_0^R \frac{r^{1-\gamma}}{\Gamma(2-\gamma)}\,dr=\frac{R^{2-\gamma}}{\Gamma(3-\gamma)}\,.
\]
In particular, for any $0<a<1$ we have
\[
    \frac{1-a^{2-\gamma}}{\Gamma(3-\gamma)}=\lim_{x\to\infty 
    (x \to 0, \text{ resp.})  
     \atop t\frac{\phi(x^{-1})^2}{
H(x^{-1})}\to 0}\frac{\int_a^1 r\P(T_t\geq xr)\,dr}{tH(x^{-1})}\geq \limsup_{x\to\infty 
    (x \to 0, \text{ resp.})  
     \atop t\frac{\phi(x^{-1})^2}{
H(x^{-1})}\to 0}\frac{(1-a^2)\P(T_t\geq x)}{2tH(x^{-1})}\,,
\]
i.e., 
\[
    \limsup_{x\to\infty 
    (x \to 0, \text{ resp.})  
     \atop t\frac{\phi(x^{-1})^2}{
H(x^{-1})}\to 0}\frac{\P(T_t\geq x)}{tH(x^{-1})}\leq \frac{2}{\Gamma(3-\gamma)}\frac{1-a^{2-\gamma}}{1-a^2}\,.
\]
Letting $a$ go to $1$ we obtain
\begin{equation}\label{eq:limsup}
\limsup_{x\to\infty 
    (x \to 0, \text{ resp.})  
     \atop t\frac{\phi(x^{-1})^2}{
H(x^{-1})}\to 0}\frac{\P(T_t\geq x)}{tH(x^{-1})}\leq \frac{2-\gamma}{\Gamma(3-\gamma)}=\frac{1}{\Gamma(2-\gamma)}\,.
\end{equation}
Similarly, for $b>1$ we obtain 
\[
    \frac{b^{2-\gamma}-1}{\Gamma(3-\gamma)}=\lim_{x\to\infty 
    (x \to 0, \text{ resp.})  
     \atop t\frac{\phi(x^{-1})^2}{
H(x^{-1})}\to 0}\frac{\int_1^b r\P(T_t\geq xr)\,dr}{tH(x^{-1})}\leq \liminf_{x\to\infty 
    (x \to 0, \text{ resp.})  
     \atop t\frac{\phi(x^{-1})^2}{
H(x^{-1})}\to 0}\frac{(b^2-1)\P(T_t\geq x)}{2tH(x^{-1})}
\]
and, by letting $b$ go to $1$, we obtain
\begin{equation}\label{eq:liminf}
\liminf_{x\to\infty 
    (x \to 0, \text{ resp.})  
     \atop t\frac{\phi(x^{-1})^2}{
H(x^{-1})}\to 0}\frac{\P(T_t\geq x)}{tH(x^{-1})}\geq \frac{1}{\Gamma(2-\gamma)}\,.
\end{equation}
The result follows now from \eqref{eq:limsup} and \eqref{eq:liminf}\,.
\qed

\begin{Rem}
It is instructive to see why previous proof does not always work if we start with \eqref{eq:tail_lt} instead of (\ref{eq:tail_lt_der}). Assume that $\phi$ varies regularly at $0$ with index $\gamma=1$\,. Using the same idea as in the proof of Theorem \ref{thm:thm_cont} we obtain 
\begin{align*}
\lim_{\tiny\begin{array}{c} x\to\infty  \\ t\phi(x^{-1})\to 0\end{array}}\frac{\int_0^\infty e^{-\lambda r} \P(T_t\geq xr)\,dr}{t\phi(x^{-1})}=\lambda^{\gamma-1}=1=\int_0^\infty e^{-\lambda r}\delta_0(dr)\,,
\end{align*}
where $\delta_0$ is the point mass at $0$, i.e.
\[
    \delta_0(A)=\begin{cases}
1\,, & 0\in A\\
0\,, & 0\not\in A\,.
\end{cases}
\]
By continuity theorem, Theorem \ref{thm:thm_cont}, we obtain 
\[
\lim_{\tiny\begin{array}{c} x\to\infty  \\ t\phi(x^{-1})\to 0\end{array}}\frac{\int_0^R\P(T_t\geq xr)\,dr}{t\phi(x^{-1})}=\int_0^R \delta_0(dr)=1\,,
\]
implying
\begin{align*}
\limsup_{\tiny\begin{array}{c} x\to\infty  \\ t\phi(x^{-1})\to 0\end{array}}\frac{\P(T_t\geq x)}{2t\phi(x^{-1})}&\leq \lim_{\tiny\begin{array}{c} x\to\infty  \\ t\phi(x^{-1})\to 0\end{array}}\frac{\int_{1/2}^1\P(T_t\geq xr)\,dr}{t\phi(x^{-1})}=1-1=0\,.
\end{align*}
Hence, in this case tail decays faster than $t\phi(x^{-1})$ as $x\to \infty$ and $t\phi(x^{-1})\to 0$\,.

\end{Rem}

\proof[Proof of Corollary \ref{cor:asympt}]
Since $1-\gamma>0$, by \eqref{eq:H_phi}, it follows that $H$ varies regularly at $0$ (at $\infty$, respectively) with index $\gamma$. Furthermore, 
\[
    t\frac{\phi(r^{-1})^2}{H(r^{-1})}\to 0 \quad \text{ if and only if } \quad t\phi(r^{-1})\to 0 \quad \text{ when }\quad r\to \infty
\]
and so, by Theorem \ref{thm:subord_tail_estimate}, 
\[
    \frac{1}{\Gamma(2-\gamma)}=\lim_{\substack{r\to\infty  (r \to 0, \text{ resp.})  \\ t\phi(r^{-1})\to 0}}\frac{\P(T_t\geq r)}{t\phi(r^{-1})}\cdot \frac{\phi(r^{-1})}{H(r^{-1})}=\frac{1}{1-\gamma}\lim_{\substack{r\to\infty (r \to 0, \text{ resp.})\\ t\phi(r^{-1})\to 0}}\frac{\P(T_t\geq r)}{t\phi(r^{-1})}\,
\]
yielding the desired estimate. 
\qed

The following proposition explains the decay of the tail of the L\' evy measure of the subordinator $T$\,.
\begin{Prop}\label{prop:levy-tail}
Assume that $H(\lambda)=\phi(\lambda)-\lambda\phi'(\lambda)$ varies regularly at $0$ (at $\infty$, respectively) with index $\gamma\in [0,1]$. Then
     \[
   \mu(r,\infty)\sim \frac{1}{\Gamma(2-\gamma)}H(r^{-1}), 
   \quad  r\to\infty \,\, (r \to 0, \text{ resp.}) \]
\end{Prop}
\proof
It is easy to show that 
\[
 b+   \int_0^\infty e^{-\lambda r}\mu(r,\infty)\,dr=\frac{\phi(\lambda)}{\lambda}\,.
\] 
Taking the derivative we obtain 
\[
    \int_0^\infty e^{-\lambda r}r\mu(r,\infty)\,dr=\frac{H(\lambda)}{\lambda^2}\,.
\]
Since
\begin{align*}
    \lambda^{-2}&=\lim_{x\to\infty \atop (x \to 0 \text{ resp.})}\int_0^\infty e^{-\lambda x^{-1}r}\frac{r\mu(r,\infty)}{x^2H(\lambda x^{-1})}\,dr\\
    &=\lim_{x\to\infty \atop (x \to 0 \text{ resp.})}\int_0^\infty e^{-\lambda r}\frac{r\mu(rx,\infty)}{H(x^{-1})}\,dr\frac{H(x^{-1})}{H(\lambda x^{-1})}\,,
\end{align*}
we can use regular varation of $H$ to get
\[
    \lim_{x\to\infty \atop (x \to 0 \text{ resp.})}\int_0^\infty e^{-\lambda r}\frac{r\mu(rx,\infty)}{H(x^{-1})}\,dr=\lambda^{\gamma-2}\,.
\]
Note that the last display is very similar to \eqref{eq:lt-tail-distr}\,. We now repeat the argument after \eqref{eq:lt-tail-distr} in the proof of Theorem \ref{thm:subord_tail_estimate}  to obtain the claim. 
\qed

\section{Asymptotic of the tail probability of isotropic L\' evy processes}

Throughout this section, we assume that $X=(X_t)_{t\geq 0}$ is an isotropic L\' evy process in $\R^d\,(d\geq 1)$ with the characteristic exponent $\xi \to \psi(|\xi|)$. The following is known and, in fact true for any
 negative definite function
 (see \cite[Lemma 1]{G}).

\begin{Lem}\label{L:1.1}
For every  $t>0$ and $\lambda \geq 1$,
$$
{\psi^*(\lambda t)} \le 2(1+\lambda^2){\psi^*(t)}
$$
where $\psi^*(r):= \sup_{|z| \le r} \psi(|z|)$.
\end{Lem}

We start with 
\begin{equation}\label{eq:gaussian}
e^{-|x|^2}=(4\pi)^{-d/2}\int_{\R^d}e^{i\xi\cdot x}e^{-\frac{|\xi|^2}{4}}\,d\xi
\end{equation}
to obtain the Laplace transform of the measure $\P(|X_t|\geq \sqrt{r})\,dr$ (see also \cite{CGT}); 
\begin{align}\nonumber
\int_0^\infty e^{-\lambda r}\P(|X_t|\geq \sqrt{r})\,dr&=\E\int_0^{|X_t|^2}e^{-\lambda r}\,dr=(4\pi)^{-d/2}\int_{\R^d}\frac{1-\E[e^{i\sqrt{\lambda}\xi\cdot X_t}]}{\lambda}e^{-\frac{|\xi|^2}{4}}\,d\xi\label{eq:lt-tail-iso}\\&=
(4\pi)^{-d/2}\int_{\R^d}\frac{1-e^{-t\psi(|\xi|\sqrt{\lambda})}}{\lambda}e^{-\frac{|\xi|^2}{4}}\,d\xi
\\&=\frac{2^{1-d}}{\Gamma(d/2)}\int_0^\infty \frac{1-e^{-t\psi(r\sqrt{\lambda})}}{\lambda}e^{-\frac{r^2}{4}}r^{d-1}\,dr\nonumber\\
&=\frac{2^{1-d}}{\Gamma(d/2)}\int_0^\infty \lambda^{-d/2-1}e^{-\frac{r^2}{4\lambda}}(1-e^{-t\psi(r)})r^{d-1}\,dr\label{eq:lt-tail-iso-polar}
\end{align}
for all $\lambda>0$\,.

One of the ideas of proving the asymptotic behaviors of the tail probability and the transition density is using the following function $\phi$, which is {\it nicer} than $\psi$.
Define the auxiliary function $\phi:(0,\infty)\rightarrow (0,\infty)$ by
\begin{equation}\label{eq:phi}
    \phi(\lambda)=(4\pi)^{-d/2}\int_{\R^d}e^{-\frac{|\xi|^2}{4}}\psi(|\xi|\sqrt{\lambda})\,d\xi\,,\quad \lambda>0\,.
\end{equation}
The following forms of $\phi$ are useful  in evaluation of its derivative (since then we do not have to take derivative of $\psi$);
\begin{equation}\label{eq:phi-polar}
    \phi(\lambda)=\int_{\R^d}(4\pi \lambda)^{-d/2}e^{-\frac{|\xi|^2}{4\lambda}}\psi(|\xi|)\,d\xi=\frac{2^{1-d}}{\Gamma(d/2)}\int_0^\infty \lambda^{-d/2}e^{-\frac{r^2}{4\lambda}}\psi(r)r^{d-1}\,dr\,.
\end{equation}
\begin{Rem}\label{rem:sbm}
The function $\phi$ is, in fact, a Bernstein function and  the Laplace exponent of a subordinator.
Moreover, if $a=0$ in \eqref{e:psirep} then $\phi$ is  the Laplace exponent of a subordinator 
 without drift. 
 
 Indeed, by using \eqref{e:psirep}, \eqref{eq:gaussian} and polar coordinates, we obtain
\begin{align*}
\phi(\lambda)&=a\int_{\R^d}(4\pi \lambda)^{-d/2}e^{-\frac{|\xi|^2}{4\lambda}}|\xi|^2d\xi+\int_{\R^d}\int_{\R^d}(4\pi \lambda)^{-d/2}e^{-\frac{|\xi|^2}{4\lambda}}(1-e^{i\xi\cdot x})\,d\xi\,\nu(dx)\\
&=a\sum_{i=1}^d 
\int_{\R^d}(4\pi \lambda)^{-d/2}e^{-\frac{|\xi|^2}{4\lambda}}\xi_i^2d\xi +\int_{\R^d}(1-e^{-\lambda|x|^2})\,\nu(dx)\\
&=a\sum_{i=1}^d 
\int_{\R}(4\pi \lambda)^{-1/2}e^{-\frac{\xi_i^2}{4\lambda}}\xi_i^2d\xi +\int_0^\infty (1-e^{-\lambda r^2})\int_{\partial B(0,r)}\nu(dy)\,dr\\
&=2a \lambda d + 2^{-1}
\int_0^\infty (1-e^{-\lambda r})r^{-1/2}\int_{\partial B(0,\sqrt{r})}\nu(dy)\,dr\,.
\end{align*}
Denoting $\mu(dr)=\frac{1}{2\sqrt r}\int_{\partial B(0,\sqrt{r})}\nu(dy)\,dr$\,, we obtain that $\mu$ is the L\' evy measure of a subordinator, since
\begin{align*}
    \int_0^\infty (1\wedge r)\mu(dr)&=
    \int_0^\infty (1\wedge r)\frac{1}{2\sqrt r}\int_{\partial B(0,\sqrt{r})}\nu(dy)\,dr
    \\&=
    \int_0^\infty (1\wedge r^2)\int_{\partial B(0,r)}\nu(dy)\,dr=\int_{\R^d}(1\wedge |x|^2)\nu(dx)<\infty\,.
\end{align*}
\end{Rem}

As in the subordinate Brownian motion case, set
\[
    H(\lambda)=\phi(\lambda)-\lambda\phi'(\lambda),\quad \lambda>0\,.
\]
By Remark \ref{rem:sbm}, 
$H$ is non-negative and non-decreasing.

 Using \eqref{eq:phi-polar}, we obtain $H$ in terms of $\psi$ as
\begin{align}
    H(\lambda)&=-\lambda^2\left(\frac{\phi(\lambda)}{\lambda}\right)'=\nonumber
    \\&=\frac{2^{1-d}}{\Gamma(\frac{d}{2})}\int_0^\infty \lambda^{-d/2}e^{-\frac{r^2}{4\lambda}}\left(-\frac{r^2}{4\lambda}+\frac{d+2}{2}\right)\psi(r)r^{d-1}\,dr\,.\label{eq:lt10}
\end{align}

We now discuss  compatibilities of 
functions  $\phi(\lambda)$,  $H(\lambda)$ and  $\psi(\sqrt{\lambda})$. 
\begin{Prop}\label{prop:hphi-regvar}
 Assume 
  that $\psi$ varies regularly at $0$ (at $\infty$, respectively) with index $\alpha\in [0,2]$. Then
  
  \noindent
  (1)
$$
 \phi(\lambda)    \sim 
 2^{\alpha}\frac{\Gamma\left(\frac{d+\alpha}{2}\right)}{\Gamma\left(\frac{d}{2}\right)}\psi(\sqrt{\lambda}), \quad \lambda \to 0 \,\,(\lambda \to \infty, \,  resp.), 
   $$
   
    \noindent
   (2)
     $$
\lim_{ \lambda \to 0 \,\,(\lambda \to \infty, \text{ resp.})}   \frac{H(\lambda)}{\psi(\sqrt{\lambda})}= 2^{\alpha-1}(2-\alpha)\frac{\Gamma\left(\frac{d+\alpha}{2}\right)}{\Gamma\left(\frac{d}{2}\right)}
 $$
  where   $\phi$ is defined in  \eqref{eq:phi-polar}. 
  Consequently, $\phi$ varies regularly at $0$ (at $\infty$, respectively) with index $\alpha/2$ and, if $\alpha<2$, then  $H$ varies regularly at $0$ (at $\infty$, respectively) with index $\alpha/2$\,.
\end{Prop}
\proof

(1) 
By \eqref{eq:phi-polar} we have 
\begin{align*}
\frac{\phi(\lambda)}{\psi(\sqrt{\lambda})}&=\frac{1}{2^{d-1}\Gamma(\frac{d}{2})}
   \left(\int_0^{1/\sqrt{\lambda}}  e^{-\frac{r^2}{4}}r^{d-1}
   \frac{\psi(r\sqrt{\lambda})}{\psi(\sqrt{\lambda})}
   \,dr+\int_{1/\sqrt{\lambda}}^\infty  e^{-\frac{r^2}{4}}r^{d-1}
   \frac{\psi(r\sqrt{\lambda})}{\psi(\sqrt{\lambda})}
   \,dr\right)
  \\
    &=  I(\lambda) +\frac{1}{2^{d-1}\Gamma(\frac{d}{2})}
 II(\lambda)
   \,. 
  \end{align*}

We first assume 
  that $\psi$ varies regularly at $0$  with index $\alpha\in [0,2]$.
 By Potter's theorem (see \cite[Theorem 1.5.6]{BGT}) there exists a constant $c_1>0$ such that 
\begin{align}
\label{e:ptu}
 \frac{ \psi(r\sqrt{\lambda})}{ \psi(\sqrt{\lambda})}\leq c_1(r^3 \vee r^{-1/2}),\quad \text{ for }\,\,1>r\sqrt{\lambda}>0, \, \, 1>\lambda>0.
\end{align}
So, by the dominated convergence theorem we have
  \begin{align*}
   \lim_{\lambda\downarrow 0}I(\lambda)&
   =\frac{1}{2^{d-1}\Gamma(\frac{d}{2})}\int_0^\infty e^{-\frac{r^2}{4}}r^{\alpha+d-1}\,dr
   \\&=\frac{1}{2^{d-1}\Gamma(\frac{d}{2})}\int_0^\infty e^{-r}2^{\alpha+d-1}r^{\frac{d+\alpha}{2}-1}\,dr=2^{\alpha}\frac{\Gamma\left(\frac{d+\alpha}{2}\right)}{\Gamma\left(\frac{d}{2}\right)}\,. 
  \end{align*}
  
 By Lemma \ref{L:1.1} we have
$\int_1^\infty e^{-\frac{r^2}{8}} \psi(r) r^{-5}\,dr =:c_2 < \infty.$
Thus $\lambda \in (0, 1]$
\begin{align}\label{e:fdaafd}
& \lambda^{-d/2}\int_{1}^\infty e^{-\frac{r^2}{4\lambda}}r^{d-1}
 \psi(r)
   \,dr= \lambda^2    \int_1^\infty \left(\frac{r}{\sqrt{\lambda}}\right)^{d+4}e^{-\frac{r^2}{4\lambda}} \psi(r)  r^{-5} dr \nonumber \\
  & \le \lambda^2    \int_1^\infty \left(\frac{r}{\sqrt{\lambda}}\right)^{d+4}e^{-\frac{r^2}{4\lambda}} \psi(r)  r^{-5} dr  \nonumber \\
&\le  \lambda^2  \big( \sup_{y \ge 1} y^{d+4}e^{-y/8}\big)\int_1^\infty e^{-\frac{r^2}{8}} \psi(r) r^{-5}\,dr
=:c_3 \lambda^2 .
 \end{align}
  Since $\lambda\mapsto \psi(\sqrt{\lambda})$ varies regularly 
with index  $\alpha/2\in [0,1]$, after change of variable, we have that 
\begin{align}\label{e:fdaafd1}
 \limsup_{\lambda\downarrow 0}II(\lambda)=
 \limsup_{\lambda\downarrow 0} \frac{\lambda^{-d/2}}{\psi(\sqrt{\lambda})}\int_{1}^\infty e^{-\frac{r^2}{4\lambda}}r^{d-1}
 \psi(r)
   \,dr 
  \leq  c_3
 \limsup_{\lambda\downarrow 0}\frac{\lambda^2}{\psi(\sqrt{\lambda})}=0.
 \end{align}
Therefore,  
\begin{align*}
   \lim_{\lambda\downarrow 0}\frac{\phi(\lambda)}{\psi(\sqrt{\lambda})}&=  \frac{1}{2^d\Gamma(\frac{d}{2})} \lim_{\lambda\downarrow 0}
      I(\lambda)=2^{\alpha}\frac{\Gamma\left(\frac{d+\alpha}{2}\right)}{\Gamma\left(\frac{d}{2}\right)}.
     \end{align*}
     
     We now assume 
  that $\psi$ varies regularly at $\infty$  with index $\alpha\in [0,2]$.
  In this case, through analogous argument
  \begin{align*}
   \lim_{\lambda\uparrow \infty}II(\lambda)=2^{\alpha}\frac{\Gamma\left(\frac{d+\alpha}{2}\right)}{\Gamma\left(\frac{d}{2}\right)}\,. 
  \end{align*} 
Moreover, by change of variable 
     \begin{align*} I(\lambda)= \frac{1}{\lambda^{d/2}\psi(\sqrt{\lambda})}\int_{0}^1 e^{-\frac{r^2}{4\lambda}}r^{d-1}
 \psi(r)
   \,dr  \le \frac{1}{\lambda^{d/2}\psi(\sqrt{\lambda})}\int_{0}^1r^{d-1}
 \psi(r)
   \,dr,
  \end{align*} 
  which goes to zero as $\lambda\uparrow \infty.$
 
 \indent 
  (2) As we observed  in (1), the proof will be analogous and simpler when $\psi$ varies regularly at $\infty$  with index $\alpha\in [0,2]$. 
Thus  we only provide the proof for the case that $\psi$ varies regularly at $0$  with index $\alpha\in [0,2]$. 

  Using \eqref{eq:lt10} and change of variable, we rewrite as
    \begin{align*}
  &\frac{H(\lambda)}{\psi(\sqrt{\lambda})}= \frac{2^{1-d}}{\Gamma(\frac{d}{2})}\int_0^\infty \lambda^{-d/2}e^{-\frac{r^2}{4\lambda}}\left(\frac{d+2}{2}-\frac{r^2}{4\lambda}\right) r^{d-1}\frac{\psi(r)}{\psi(\sqrt{\lambda})}dr\\
  &= \frac{2^{1-d}}{\Gamma(\frac{d}{2})}\int_0^1 \lambda^{-d/2}e^{-\frac{r^2}{4\lambda}}\left(\frac{d+2}{2}-\frac{r^2}{4\lambda}\right) r^{d-1}\frac{\psi(r)}{\psi(\sqrt{\lambda})}dr\\
  &\quad + \frac{2^{1-d}}{\Gamma(\frac{d}{2})}\int_1^\infty \lambda^{-d/2}e^{-\frac{r^2}{4\lambda}}\left(\frac{d+2}{2}-\frac{r^2}{4\lambda}\right) r^{d-1}\frac{\psi(r)}{\psi(\sqrt{\lambda})}dr\\
   &=\frac{1}{\Gamma(\frac{d}{2})}\int_0^{(4\lambda)^{-1}} e^{-r}r^{\frac{d}{2}-1}\left(\frac{d+2}{2}-r\right)\frac{\psi(2\sqrt{r}\sqrt{\lambda})}{\psi(\sqrt{\lambda})}\,dr\\
&\quad + \frac{2^{1-d}}{\Gamma(\frac{d}{2})}\frac{1}{ \psi(\sqrt{\lambda}) }\int_1^\infty \lambda^{-d/2}e^{-\frac{r^2}{4\lambda}}\left(\frac{d+2}{2}-\frac{r^2}{4\lambda}\right) r^{d-1}\psi(r)dr \\&=I(\lambda)+\frac{2^{1-d}}{\Gamma(\frac{d}{2}) } II(\lambda)\, .
  \end{align*}
  Similar to (1), by Potter's theorem as \eqref{e:ptu} we can use the dominated convergence theorem and get 
   \begin{align*}
   \lim_{\lambda\downarrow 0}  I(\lambda) &=\frac{1}{\Gamma(\frac{d}{2})}\int_0^\infty e^{-r}r^{\frac{d}{2}-1}\left(\frac{d+2}{2}-r\right)2^{\alpha}r^{\frac{\alpha}{2}}\,dr\\
   &=\frac{2^\alpha}{\Gamma(\frac{d}{2})}\left(\frac{d+2}{2}\Gamma\left(\frac{d+\alpha}{2}\right)-\Gamma\left(\frac{d+\alpha}{2}+1\right)\right)\\
   &=\frac{2^\alpha}{\Gamma\left(\frac{d}{2}\right)}\left(\frac{d+2}{2}\Gamma\left(\frac{d+\alpha}{2}\right)-\frac{d+\alpha}{2}\Gamma\left(\frac{d+\alpha}{2}\right)\right)\\
   &=2^{\alpha-1}(2-\alpha)\frac{\Gamma\left(\frac{d+\alpha}{2}\right)}{\Gamma\left(\frac{d}{2}\right)}.
  \end{align*}
  Similar as \eqref{e:fdaafd} and \eqref{e:fdaafd1}, for $\lambda \in (0, 1]$ we also bound  $II(\lambda)$ as
     \begin{align*}
 &  |II(\lambda)| \le 
   \frac{1}{ \psi(\sqrt{\lambda}) }\int_1^\infty \lambda^{-d/2}e^{-\frac{r^2}{4\lambda}}\left(\frac{d+2}{2}+\frac{r^2}{4\lambda}\right) r^{d-1}\psi(r)dr \\
  & \le \frac{2d+5}{4} \frac{1}{ \psi(\sqrt{\lambda}) }\int_1^\infty \lambda^{-d/2-1}e^{-\frac{r^2}{4\lambda}}r^{d+1}\psi(r)dr \\
   &= \frac{2d+5}{4} \frac{\lambda^2}{ \psi(\sqrt{\lambda}) }
   \int_1^\infty \left(\frac{r}{\sqrt{\lambda}}\right)^{d+6}e^{-\frac{r^2}{4\lambda}}r^{-5}\psi(r)dr \\
   & \le \frac{2d+5}{4} \frac{\lambda^2}{ \psi(\sqrt{\lambda}) }\big( \sup_{y \ge 1} y^{d+6}e^{-y/8}\big)
   \int_1^\infty e^{-\frac{r^2}{8}}r^{-5}\psi(r)dr, 
     \end{align*}
     which goes to zero as $\lambda \to 0$. 
     
     Therefore, 
      \begin{align*}
  \lim_{\lambda\downarrow 0}\frac{H(\lambda)}{\psi(\sqrt{\lambda})}=2^{\alpha-1}(2-\alpha)\frac{\Gamma\left(\frac{d+\alpha}{2}\right)}{\Gamma\left(\frac{d}{2}\right)}.
  \end{align*}

  \qed

   It follows from the previous result that, in the case $\alpha=2$, functions $H(\lambda)$ and $\psi(\sqrt{\lambda})$ are not comparable. We clarify this situation in the following result. 
  \begin{Prop}\label{prop:Hregvar2}
  Let $g(\lambda)=\psi(\lambda)-\frac{1}{2}\lambda\psi'(\lambda)$. If $g$ varies reularly at $0$ (at $\infty$, respectively) with index $\alpha \ge 0$, then
  $$
H(\lambda)=\phi(\lambda)-\lambda\phi'(\lambda)  \sim 
 \frac{2^\alpha\Gamma\left(\frac{d+\alpha}{2}\right)}{\Gamma\left(\frac{d}{2}\right)}g(\sqrt{\lambda}), \quad \lambda \to 0 \,\,(\lambda \to \infty, \text{resp.}).
   $$

  \end{Prop}
\proof
Since the proof for the case that $g$ varies regularly at $\infty$  will be analogous and simpler,   we only provide the proof for the case that $g$ varies regularly at $0$  with index $\alpha \ge0$. 
%
%\noindent
%(1) Suppose  $g$ varies regularly at $0$  with index $\alpha \ge0$. 

Note that, by the change of variable we have  
    \begin{align}
\phi(\lambda) &= \frac{2^{1-d}}{\Gamma(\frac{d}{2})}  \left(\int_0^{1/\sqrt{\lambda}}e^{-r^2/4}r^{d-1}{\psi(r\sqrt{\lambda})}dr+ 
 \int_{1/\sqrt{\lambda}}^\infty e^{-r^2/4}r^{d-1}{\psi(r\sqrt{\lambda})}dr \right)\nonumber\\
&=  \frac{2^{1-d}}{\Gamma(\frac{d}{2})}  \left(\int_0^{1/\sqrt{\lambda}}e^{-r^2/4}r^{d-1}{\psi(r\sqrt{\lambda})}dr+ 
  \lambda^{-d/2} \int_1^\infty r^{d-1}e^{-\frac{r^2}{4\lambda}}\psi(r) dr
\right).
 \label{eq:lt100}
\end{align}

   By Potter's theorem (see \cite[Theorem 1.5.6]{BGT}) there exists a constant $c_0>0$ such that 
    \[
     {g(r\sqrt{\lambda})}\leq c_0{g(1)}(r\sqrt{\lambda})^{-1/2},\quad \text{ for }\,\,1>r\sqrt{\lambda}>0, \, \, 
    \]
Thus we can differentiate under the  integral sign and get 
\begin{align}
&\left( \int_0^{1/\sqrt{\lambda}}e^{-r^2/4}r^{d-1}\lambda^{-1} {\psi(r\sqrt{\lambda})}dr \right)' \nonumber\\
=&      e^{-1/(4\lambda)} \lambda^{-(d-1)/2} \lambda^{-1}\psi(1) (-\frac{1}{2} \lambda^{-3/2} )+   \int_0^{1/\sqrt{\lambda}}e^{-r^2/4}r^{d-1} 
\frac{
\partial} {\partial \lambda} \left(\frac{\psi(r\sqrt{\lambda})} {\lambda} \right)dr     \nonumber\\    
=&-\frac{1}{2}e^{-1/(4\lambda)} \lambda^{-2-d/2} \psi(1) -   \int_0^{1/\sqrt{\lambda}}e^{-r^2/4}r^{d-1} \frac{g(r\sqrt{\lambda})} {\lambda^2}dr.
 \label{eq:lt1000} 
\end{align}

    From \eqref{eq:lt100} and \eqref{eq:lt1000}  we obtain that 
    \begin{align}
    &H(\lambda)=-\lambda^2\left(\frac{\phi(\lambda)}{\lambda}\right)'\nonumber
    \\&=- \frac{2^{1-d}}{\Gamma(\frac{d}{2})} \lambda^2\left(   
   \int_0^{1/\sqrt{\lambda}}e^{-r^2/4}r^{d-1}\lambda^{-1}{\psi(r\sqrt{\lambda})}dr+ 
  \lambda^{-1-d/2} \int_1^\infty r^{d-1}e^{-\frac{r^2}{4\lambda}}\psi(r) dr
 \right)'\nonumber\\
   &=- \frac{2^{1-d}}{\Gamma(\frac{d}{2})} \lambda^2\left(    
   -\frac{1}{2}e^{-1/(4\lambda)} \lambda^{-2-d/2} \psi(1) -   \int_0^{1/\sqrt{\lambda}}e^{-r^2/4}r^{d-1} \frac{g(r\sqrt{\lambda})} {\lambda^2}dr \right.\nonumber\\
  &\quad \left. +
\int_1^\infty \lambda^{-2-d/2}e^{-\frac{r^2}{4\lambda}}\left(\frac{r^2}{4\lambda}-\frac{d+2}{2}\right)\psi(r)r^{d-1}\,dr \right)\nonumber\\
   &= \frac{2^{1-d}}{\Gamma(\frac{d}{2})} \left( \int_0^{1/\sqrt{\lambda}}e^{-r^2/4}r^{d-1} {g(r\sqrt{\lambda})} dr    
 +\frac{1}{2}e^{-1/(4\lambda)} \lambda^{-d/2} \psi(1)\right.\nonumber\\
  &\quad \left. +
\int_1^\infty \lambda^{-d/2}e^{-\frac{r^2}{4\lambda}}\left(-\frac{r^2}{4\lambda}+\frac{d+2}{2}\right)\psi(r)r^{d-1}\,dr \right).\label{eq:lt10}
\end{align}
   
    To evaluate the limit we use the above form so that 
    \begin{align*}
       & \frac{H(\lambda)}{g(\sqrt{\lambda})}=\frac{2^{1-d}}{\Gamma(\frac{d}{2})}\int_0^{1/\sqrt{\lambda}}e^{-r^2/4}r^{d-1}\frac{g(r\sqrt{\lambda})}{g(\sqrt{\lambda})}\,dr\\
        &+\frac{2^{1-d}}{\Gamma(\frac{d}{2})g(\sqrt{\lambda})}
       \left( \frac{1}{2}e^{-1/(4\lambda)} \lambda^{-d/2} \psi(1)+\int_1^\infty \lambda^{-d/2}e^{-\frac{r^2}{4\lambda}}\left(-\frac{r^2}{4\lambda}+\frac{d+2}{2}\right)\psi(r)r^{d-1}\,dr \right)
        \\&=I_1(\lambda)+\frac{2^{1-d}}{\Gamma(\frac{d}{2})g(\sqrt{\lambda})}I_2(\lambda).
    \end{align*}
    By Potter's theorem (see \cite[Theorem 1.5.6]{BGT}) there exists a constant $c_1>0$ such that 
    \[
        \frac{g(r\sqrt{\lambda})}{g(\sqrt{\lambda})}\leq c_1(r^{\alpha+1} \vee r^{-1/2}),\quad \text{ for }\,\,1>r\sqrt{\lambda}>0, \, \, 1>\lambda>0
    \]
    and so we may use dominated convergence theorem in the first integral to obtain
    \[
        \lim_{\lambda\downarrow 0} I_1(\lambda)= \frac{2^{1-d}}{\Gamma(\frac{d}{2})}\int_0^\infty e^{-r^2/4}r^{d+\alpha-1}\,dr=\frac{2^\alpha}{\Gamma(\frac{d}{2})}\int_0^\infty e^{-r}r^{\frac{d+\alpha}{2}-1}\,dr=\frac{2^\alpha\Gamma\left(\frac{d+\alpha}{2}\right)}{\Gamma\left(\frac{d}{2}\right)}.
    \]

 On the other hand, 
      we can rewrite $I_2$ as 
   \begin{align*}
    I_2(\lambda)    =&(1+d/2)  \lambda^{-d/2}  \int_1^\infty r^{d-1}e^{-\frac{r^2}{4\lambda}} \psi(r)  dr    \\
      &+\frac{1}{2}  \lambda^{-d/2}e^{-\frac{1}{4\lambda}}\psi(1) -  \frac{1}{4} \lambda^{-1-d/2}   \int_1^\infty r^{d+1}e^{-\frac{r^2}{4\lambda}}  
     \psi(r)dr\\
    =&\lambda^{(1+\alpha)/2}2\left(    (1+d/2)     \int_1^\infty \left(\frac{r}{\sqrt{\lambda}}\right)^{d+1+\alpha}e^{-\frac{r^2}{4\lambda}} \psi(r)  r^{-2-\alpha} dr    \right. \\
      & \left. + \frac{1}{2}  \lambda^{-(d+1+\alpha)/2}e^{-\frac{1}{4\lambda}}\psi(1) -  \frac{1}{4}   \int_1^\infty \left(\frac{r}{\sqrt{\lambda}}\right)^{d+3+\alpha}e^{-\frac{r^2}{4\lambda}}  
     \psi(r) r^{-2-\alpha}dr
\right)
              \end{align*}
              Thus for $\lambda \in (0,1]$
    \begin{align*}
 &   |I_2(\lambda)| \le c_2\lambda^{(1+\alpha)/2} \left(  \lambda^{-(d+1+\alpha)/2}e^{-\frac{1}{4\lambda}}+ \int_1^\infty \left(\frac{r}{\sqrt{\lambda}}\right)^{d+3+\alpha}e^{-\frac{r^2}{4\lambda}}\psi(r) r^{-2-\alpha} dr \right)\\
    &\leq c_2\lambda^{(1+\alpha)/2}  \left( \big( \sup_{y \ge 1} y^{(d+1+\alpha)/2}e^{-y/4}\big) +\big( \sup_{y \ge 1} y^{d+3+\alpha}e^{-y/8}\big)\int_1^\infty e^{-\frac{r^2}{8}} \psi(r) r^{-2-\alpha}\,dr \right)\\
    &\leq c_3\lambda^{(1+\alpha)/2}, 
    \end{align*}
       where in the last inequality  we have used $\psi(r)\leq 2(1+r^2)\sup_{s \le 1}\psi(s)$ for all $r>0$ by   Lemma \ref{L:1.1}.
 Since $\lambda\mapsto g(\sqrt{\lambda})$ varies regularly 
with index  $\alpha/2$ we now conclude that 
 \[
\limsup_{\lambda\downarrow 0}\frac{|I_2(\lambda)|}{g(\sqrt{\lambda})}\leq c_3\limsup_{\lambda\downarrow 0}\frac{\lambda^{(1+\alpha)/2}}{g(\sqrt{\lambda})}=0.
 \]

\qed

\begin{Lem}\label{l:main}
 Assume that $H(\lambda)=\phi(\lambda)-\lambda\phi'(\lambda)$ varies regularly at $0$ (at $\infty$, respectively) with index 
 %$\gamma\in [0,1]$.
 $\gamma\in [0,2)$ and $\psi$ varies regularly at $0$ (at $\infty$, respectively)  with index $\gamma_1 \in [0,2]$.
   Then  for all $R>0$ we have 
\begin{align*}
 \int_0^Rr\P(|X_t|\geq \sqrt{r}\sqrt{s})\,dr\sim \frac{R^{2-\gamma}}{\Gamma(3-\gamma)} tH(s^{-1}), 
  & \quad  s\to\infty \,\, (s \to 0, \text{ resp.})\\
  \text{ and }\quad & t\frac{\psi(\sqrt{s}^{-1})^2}{H(s^{-1})} \to 0\,.
\end{align*}
\end{Lem}

\proof
We provide the proof for the case  that $H$ varies regularly at $0$ with index 
 $\gamma\in [0,2)$. The proof of the other case is similar. 

Taking derivative in \eqref{eq:lt-tail-iso-polar} we obtain
\begin{align*}
 \int_0^\infty e^{-\lambda r}r&\P(|X_t|\geq \sqrt{r})\,dr\\&=\frac{2^{1-d}}{\Gamma(d/2)}\int_0^\infty \lambda^{-d/2-2}e^{-\frac{r^2}{4\lambda}}\left(\frac{d+2}{2}-\frac{r^2}{4\lambda}\right)(1-e^{-t\psi(r)})r^{d-1}\,dr\,.
\end{align*}

Since  we have from 
\eqref{eq:lt10}
\begin{align*}
 \frac{H(\lambda)}{\lambda^2}=\frac{2^{1-d}}{\Gamma(d/2)}\int_0^\infty \lambda^{-d/2-2}e^{-\frac{r^2}{4\lambda}}\left(\frac{d+2}{2}-\frac{r^2}{4\lambda}\right)\psi(r)r^{d-1}\,dr, 
\end{align*}
comparing the above two displays and changing variable
we see that 
\begin{align}
 \int_0^\infty& e^{-\lambda r}r\P(|X_t|\geq \sqrt{r})\,dr=\frac{tH(\lambda)}{\lambda^2}\nonumber\\
    &-\frac{2^{1-d}}{\Gamma(d/2)}\int_0^\infty e^{-\frac{r^2}{4}}\left(\frac{r^2}{4}-\frac{d+2}{2}\right)\frac{1-e^{-t\psi(r\sqrt{\lambda})}-t\psi(r\sqrt{\lambda})}{\lambda^2}r^{d-1}\,dr\,.\label{eq:lt-tail0}
\end{align}
We replace $\lambda$ in \eqref{eq:lt-tail0}  by  $\lambda s^{-1}$ and divide both sides by $tH(s^{-1})$ and get 
\begin{align*}
 \int_0^\infty &r e^{-\lambda s^{-1} r}\frac{\P(|X_t|\geq \sqrt{r})}{t H(s^{-1})}\,dr=s^2\frac{H(\lambda s^{-1})}{\lambda^2H(s^{-1})}\nonumber\\
    &-s^2\frac{2^{1-d}}{\Gamma(d/2)}\int_0^\infty e^{-\frac{r^2}{4}}\left(\frac{r^2}{4}-\frac{d+2}{2}\right)\frac{1-e^{-t\psi(r\sqrt{\lambda}\sqrt{s}^{-1})}-t\psi(r\sqrt{\lambda}\sqrt{s}^{-1})}{\lambda^2tH(s^{-1})}r^{d-1}\,dr\,.
\end{align*}
Now  changing variable in the first integral yields
\begin{align}
& \int_0^\infty  r e^{-\lambda r}\frac{\P(|X_t|\geq \sqrt{r}\sqrt{s})}{tH(s^{-1})}\,dr=\frac{H(\lambda s^{-1})}{\lambda^2H(s^{-1})}\nonumber\\&-\frac{2^{1-d}}{\Gamma(d/2)}\int_0^\infty e^{-\frac{r^2}{4}}\left(\frac{r^2}{4}-\frac{d+2}{2}\right)\frac{1-e^{-t\psi(r\sqrt{\lambda}\sqrt{s}^{-1})}-t\psi(r\sqrt{\lambda}\sqrt{s}^{-1})}{\lambda^2tH(s^{-1})}r^{d-1}\,dr\,.\label{eq:intterm1}
\end{align}
Note that, since
\begin{align*}
 \frac{|1-e^{-t\psi(r\sqrt{\lambda}\sqrt{s}^{-1})}-t\psi(r\sqrt{\lambda}\sqrt{s}^{-1})|}{tH(s^{-1})}&\leq \frac{1}{2}\frac{t\psi(r\sqrt{\lambda}\sqrt{s}^{-1})^2}{H(s^{-1})}\\&=\frac{1}{2}\left(\frac{\psi(r\sqrt{\lambda}\sqrt{s}^{-1})}{\psi(\sqrt{s}^{-1})}\right)^2\frac{t\psi(\sqrt{s}^{-1})^2}{H(s^{-1})}\,,
\end{align*}
we have that for $s \ge \lambda$
\begin{align*}
&\left|\int_0^\infty e^{-\frac{r^2}{4}}\left(\frac{r^2}{4}-\frac{d+2}{2}\right)\frac{1-e^{-t\psi(r\sqrt{\lambda}\sqrt{s}^{-1})}-t\psi(r\sqrt{\lambda}\sqrt{s}^{-1})}{\lambda^2tH(s^{-1})}r^{d-1}\,dr\right|\\
\le& \frac{1}{2}\frac{t\psi(\sqrt{s}^{-1})^2}{H(s^{-1})} \int_0^\infty e^{-\frac{r^2}{4}}\left|\frac{r^2}{4}-\frac{d+2}{2}\right|\left(\frac{\psi(r\sqrt{\lambda}\sqrt{s}^{-1})}{\psi(\sqrt{s}^{-1})}\right)^2
r^{d-1}\,dr\\
\le& \frac{1}{2}\frac{t\psi(\sqrt{s}^{-1})^2}{H(s^{-1})} \left( \int_0^{\sqrt{s}\sqrt{\lambda}^{-1}}  e^{-\frac{r^2}{4}}\left|\frac{r^2}{4}-\frac{d+2}{2}\right|\left(\frac{\psi(r\sqrt{\lambda}\sqrt{s}^{-1})}{\psi(\sqrt{s}^{-1})}\right)^2 r^{d-1}\,dr \right. \\
& \left.+\frac{2d+5}{4}\int_{\sqrt{s}\sqrt{\lambda}^{-1}}^\infty e^{-\frac{r^2}{4}}\left(\frac{\psi(r\sqrt{\lambda}\sqrt{s}^{-1})}{\psi(\sqrt{s}^{-1})}\right)^2
r^{d+1}\,dr \right)
\\
=:& \frac{1}{2}\frac{t\psi(\sqrt{s}^{-1})^2}{H(s^{-1})} \left( I(s)+ \frac{2d+5}{4}II(s)\right).
\end{align*}
Since $\psi$ varies regularly at $0$ with index $\gamma_1$,  by Potter's theorem  (see \cite[Theorem 1.5.6]{BGT}) it follows that, for any $\varepsilon>0$ there exists $c=c(\varepsilon)>0$ such that  
\begin{equation}\label{eq:potter2}
 \frac{\psi(r\sqrt{\lambda}\sqrt{s}^{-1})}{\psi(\sqrt{s}^{-1})}\leq c((r\sqrt{\lambda})^{\gamma_1-\varepsilon}+(r\sqrt{\lambda})^{\gamma_1+\varepsilon})\,
 \quad \text{for } s\ge 1, 0<r \le \sqrt{s}\sqrt{\lambda}^{-1}.
\end{equation}

Hence, by \eqref{eq:potter2} we may use dominated convergence theorem to $I(s)$ and get 
$$\lim_{s \to \infty} I(s) = \lambda^{\gamma_1/2} \int_0^{\infty}  e^{-\frac{r^2}{4}}\left|\frac{r^2}{4}-\frac{d+2}{2}\right| r^{\gamma_1+d-1}\,dr <\infty.
$$
On the other hand, by the  change of variable $v=r \sqrt{\lambda} \sqrt{s}^{-1}$ we have that
\begin{align*}
II(s)&=\int_{\sqrt{s}\sqrt{\lambda}^{-1}}^\infty e^{-\frac{r^2}{4}}\left(\frac{\psi(r\sqrt{\lambda}\sqrt{s}^{-1})}{\psi(\sqrt{s}^{-1})}\right)^2
r^{d+1}\,dr\\&=\psi(\sqrt{s}^{-1})^{-2}
\int_{1}^\infty e^{-\frac{sv^2}{4\lambda}} \psi(v)^2
 (\sqrt{s}\sqrt{\lambda}^{-1})^{d+2}\, v^{d+1} dv
\end{align*}
By Lemma \ref{L:1.1}, for $s \ge \lambda$
\begin{align*}
&\int_{1}^\infty e^{-\frac{sv^2}{4\lambda}} \psi(v)^2
 (\sqrt{s}\sqrt{\lambda}^{-1})^{d+2}\, v^{d+1} dv \\
 \le & 8 \left(\sup_{a \le 1} \psi(a)^2 \right) \int_{1}^\infty e^{-\frac{sv^2}{4\lambda}}  (\sqrt{s}\sqrt{\lambda}^{-1})^{d+2}\, v^{d+5} dv\\
=& 8 \left(\sup_{a \le 1} \psi(a)^2 \right) s^{-3}\lambda^3\int_{1}^\infty e^{-\frac{sv^2}{4\lambda}}  (v\sqrt{s}\sqrt{\lambda}^{-1})^{d+8}\,v^{-3} dv\\
\le & 8 \left(\sup_{a \le 1} \psi(a)^2 \right) s^{-3}\lambda^3  \left(\sup_{b \ge 1} b^{d+8} e^{-\frac{b^2}{4}}   \right) 
\int_{1}^\infty  v^{-3} dv\\
= & 4 \left(\sup_{a \le 1} \psi(a)^2 \right) s^{-3}\lambda^3  \left(\sup_{b \ge 1} b^{d+8} e^{-\frac{b^2}{4}}   \right).
\end{align*}
Thus, 
since $a \to \psi(\sqrt a^{-1})$ varies regularly at $\infty$ with index $\gamma_1/2 \in [0,1]$, 
\begin{align*}
0\le \limsup_{s \to \infty} II(s)&\le4  \lambda^3 \left(\sup_{a \le 1} \psi(a)^2 \right) \left(\sup_{b \ge 1} b^{d+6} e^{-\frac{b^2}{8}}   \right) \left(\limsup_{s \to \infty} \frac{s^{-3/2}}{\psi(\sqrt{s}^{-1})} \right)^2 =0
\end{align*}
Therefore we conclude that the integral term in \eqref{eq:intterm1} goes to $0$ if $s\to\infty$ and $ t\frac{\psi(\sqrt{s}^{-1})^2}{H(s^{-1})} \to 0$\,.
This and regular variation of $H$ with index $\gamma$ (by Proposition \ref{prop:Hregvar2}) yield
\[
  \lim_{\substack{s\to\infty\\t\frac{\psi(\sqrt{s}^{-1})^2}{H(s^{-1})}\to 0}}\int_0^\infty e^{-\lambda r}\frac{r\P(|X_t|\geq \sqrt{r}\sqrt{s})}{tH(s^{-1})}\,dr=\lambda^{\gamma-2}=\int_0^\infty e^{-\lambda r}\frac{r^{1-\gamma}}{\Gamma(2-\gamma)}\,dr\,.
\]
By Theorem \ref{thm:thm_cont} it follows that 
\[
 \lim_{\substack{s\to\infty\\t\frac{\psi(\sqrt{s}^{-1})^2}{H(s^{-1})}\to 0}}\int_0^R \frac{r\P(|X_t|\geq \sqrt{r}\sqrt{s})}{tH(s^{-1})}\,dr=\int_0^R\frac{r^{1-\gamma}}{\Gamma(2-\gamma)}\,dr=\frac{R^{2-\gamma}}{\Gamma(3-\gamma)}
\]
for any $R>0$\,.
\qed

\begin{Prop}\label{prop:tail}
 Assume that $H(\lambda)=\phi(\lambda)-\lambda\phi'(\lambda)$ varies regularly at $0$ (at $\infty$, respectively) with index 
 $\gamma\in [0,2)$
 and $\psi$ varies regularly at $0$ (at $\infty$, respectively)  with index $\gamma_1 \in [0,2]$.
Then
 $$
\P(|X_t|\geq r) \sim\frac{1}{\Gamma(2-\gamma)}tH(r^{-2}),
 \quad  r\to\infty \,\, (r \to 0, \text{ resp.})
  \text{ and }\,  t\frac{\psi(r^{-1})^2}{H(r^{-2})} \to 0\,.
$$
\end{Prop}
\proof
By Lemma \ref{l:main},  we have that for $0<a<b$ 
    \begin{align*}
    \lim_{\substack{r\to\infty  (r \to 0, \text{ resp.})\\t\frac{\psi(\sqrt{r}^{-1})^2}{H(r^{-1})}\to 0}}\frac{\int_a^bs\P(|X_t|\geq \sqrt{rs})\,ds}{tH(r^{-1})}=\frac{b^{2-\gamma}-a^{2-\gamma}}{\Gamma(3-\gamma)},
    \end{align*}
which implies that 
\begin{align*}
\limsup_{\substack{r\to\infty  (r \to 0, \text{ resp.})\\t\frac{\psi(\sqrt{r}^{-1})^2}{H(r^{-1})}\to 0}}\frac{b^2-a^2}{2}\frac{\P(|X_t|\geq \sqrt{br})}{tH(r^{-1})}&\leq \frac{b^{2-\gamma}-a^{2-\gamma}}{\Gamma(3-\gamma)}\\
&\leq \limsup_{\substack{r\to\infty  (r \to 0, \text{ resp.})\\t\frac{\psi(\sqrt{r}^{-1})^2}{H(r^{-1})}\to 0}}\frac{b^2-a^2}{2}\frac{\P(|X_t|\geq \sqrt{ar})}{tH(r^{-1})}.
\end{align*}
By taking $b=1>a$ we get
\[
\limsup_{\substack{r\to\infty  (r \to 0, \text{ resp.})\\t\frac{\psi(\sqrt{r}^{-1})^2}{H(r^{-1})}\to 0}}\frac{\P(|X_t|\geq \sqrt{
r})}{tH(r^{-1})}\leq \lim_{a\uparrow 1}\frac{1-a^{2-\gamma}}{1-a^2}\frac{2}{\Gamma(3-\gamma)}=\frac{1}{\Gamma(2-\gamma)}.
\]
Similarly, for $a=1<b$ we obtain
\[
\liminf_{\substack{r\to\infty  (r \to 0, \text{ resp.})\\t\frac{\psi(\sqrt{r}^{-1})^2}{H(r^{-1})}\to 0}}\frac{\P(|X_t|\geq \sqrt{
r})}{tH(r^{-1})}\geq \lim_{b\downarrow 1}\frac{b^{2-\gamma}-1}{b^2-1}\frac{2}{\Gamma(3-\gamma)}=\frac{1}{\Gamma(2-\gamma)}.
\]
\qed

We now observe the following fact.
\begin{Prop}\label{prop:Hregvarn}
 Assume that $g(\lambda)=\psi(\lambda)-\frac{\lambda}2\psi'(\lambda)$ varies regularly at $0$ (at $\infty$, respectively) with index 
 $2\gamma$ with $\gamma \ge 0$. We further assume that, if $g$ varies regularly at $\infty$ with index 
 $2\gamma \in [0, 2)$, then the diffusion part is zero. 
 Then $\psi$ varies regularly at $0$ (at $\infty$, respectively) with index $2(1  \wedge \gamma)$\,.
\end{Prop}
\proof

Since $$\left(\log\frac{\psi(t)}{t^2}\right)'=\frac{\psi'(t)}{\psi(t)}-\frac{2}{t}      =-\frac{2g(t)}{t \psi(t)},$$ by integrating it we obtain
\[
 \log\frac{\psi(\lambda)}{\lambda^2}=\log\psi(1)+\int_\lambda^1\frac{2g(t)}{t \psi(t)}\,dt\,.
\]
Hence,
\begin{align}\label{e:pHrepn}
 \frac{\psi(\lambda)}{\lambda^2}=\psi(1)\exp \left(\int_\lambda^1\frac{2g(t)}{ \psi(t)} \,\frac{dt}{t} \right)\quad \text{ for}\,\,\lambda >0\,,
\end{align}

On the other hand, since $(\frac{\psi(\lambda)}{\lambda^2})'=-2\lambda^{-3}g(\lambda)$, 
we have
\begin{align}
\label{e:fplfn}
 \frac{\psi(\lambda)}{\lambda^2}-\frac{\psi(R)}{R^2}=\int_\lambda^R \frac{2g(s)}{s^3}\,ds
\end{align}
Thus for $\lambda >0$ and $\theta \ge 1$
\begin{align}
\label{e:fplf1n}
\frac{\psi(\lambda)}{g(\lambda)}=
 \frac{\psi(\lambda)/\lambda^2}{g(\lambda)/\lambda^2}=
 \int_\lambda^{\theta \lambda}\frac{2\lambda^2 g(s)}{s^3g(\lambda)}\,ds+\frac{ \psi(\theta \lambda )}{\theta^2 g(\lambda)}
 = \int_1^{\theta}\frac{2g(\lambda s)}{s^3g(\lambda)}\,ds+ \frac{ \psi(\theta \lambda )}{\theta^2 g(\lambda)}.
\end{align}

\noindent
(1) We first assume $g$ varies regularly at  $0$ with index 
 $2\gamma$ with $\gamma \in [0, 1)$. Then, by Potter's theorem (see \cite[Theorem 1.5.6]{BGT}), there exists a constant $c_1>0$ such that 
\begin{align*}
\frac{ g(\lambda s)}{s^3g(\lambda)}\leq c_1 
s^{\gamma-2}
\quad \text{ for }\,\, 0<\lambda<\lambda  s<1 \text{ and } s>1.
\end{align*}
Thus by the dominated convergence theorem,
from \eqref{e:fplf1n} we have 
\begin{align*}
\lim_{\lambda \to 0} \frac{\psi(\lambda)}{g(\lambda)}=
\lim_{\lambda \to 0}\left( \int_1^{1/\lambda}\frac{2g(\lambda s)}{s^3g(\lambda)}\,ds+ \psi(1 )\frac{\lambda^2}{g(\lambda)}\right)=
 \int_1^{\infty}2s^{2\gamma-3}\,ds=\frac{1}{1-\gamma}\end{align*}
 Therefore 
$
\lim_{\lambda \to 0} \frac{2g(\lambda)}{\psi(\lambda)}={2-2\gamma}$.
Now, using the zero version of  \cite[Theorem 1.3.1 and (1.5.2)]{BGT} we conclude from \eqref{e:pHrepn}  and the above display that $\psi(\lambda)/{\lambda^2}$ varies  regularly at $0$ with index $2\gamma-2$\,. Therefore, $\psi$ varies regularly at $0$ with index $2\gamma$\,.

\noindent
(2)
Assume that  $g$ varies regularly at  $\infty$ with index 
 $2\gamma$ with $\gamma \in [0, 1)$ and  the diffusion part  is zero. Then 
\[
 \lim\limits_{\lambda\to\infty}\frac{\psi(\lambda)}{\lambda^2}=0
\]
 Hence, by letting $R\to\infty$ in \eqref{e:fplfn}
we obtain
\begin{align}
\label{e:pll}
\frac{\psi(\lambda)}{\lambda^2}=\int_\lambda^\infty \frac{2g(s)}{s^2}\,\frac{ds}{s}\,.
\end{align}
\eqref{e:pll} and \cite[Proposition1.5.10]{BGT} imply that $\psi$ also varies regularly at $\infty$ with index 
 $2\gamma$ with $\gamma \in [0, 1)$.

\noindent
(3) Assume that  $g$ varies regularly at $0$  (at $\infty$, respectively) with index  $2\gamma \ge 2$. 
Proposition \ref{prop:Hregvar2} implies that the term $\frac{ \psi(\theta \lambda )}{\theta^2 g(\lambda)}$ in \eqref{e:fplf1n} stays positive for small $\lambda$ (large $\lambda$, respectively). Thus, from \eqref{e:fplf1n}, we have that 
 for any $\theta\geq 1$ 
\[
 \liminf_{\lambda\downarrow 0\atop 
(\lambda\uparrow \infty, \text{ resp.}) }
\frac{\psi(\lambda)}{g(\lambda)}
\ge 
2 \int_1^\theta s^{2\gamma-3}  ds=
\begin{cases} 2 \log\theta& \text{if } \gamma=1\,\\
\frac1{\gamma -1}(\theta^{2(\gamma -1)}-1)&\text{if } \gamma>1 \,,
\end{cases}
\]
where  we have used the fact
\[
 \lim_{\lambda\downarrow 0\atop 
(\lambda\uparrow \infty, \text{ resp.}) }
\frac{g(\lambda s)
}{{g(\lambda)}}=s^{2\gamma}\quad \text{ uniformly in }\,\,\,s\in [1,\theta].
\]
 (See \cite[Theorem 1.5.2]{BGT}.) 
Therefore, by letting $\theta\to\infty$, we get 
\begin{align}\label{e:fdreq1}
 \lim_{t\downarrow 0\atop 
(t\uparrow \infty, \text{ resp.}) } \frac{g(t)}{ \psi(t)} = \left(
 \lim_{t\downarrow 0\atop 
(t\uparrow \infty, \text{ resp.}) }
 \frac{{\psi(t)}}{{g(t)}} \right)^{-1}=0\,.
\end{align}
 Using \cite[Theorem 1.3.1]{BGT} we conclude from \eqref{e:pHrepn}  and \eqref{e:fdreq1} that $\frac{\psi(\lambda)}{\lambda^2}$ varies slowly at $0$ (at $\infty$, respectively). Therefore, $\psi$ varies regularly at $0$ (at $\infty$, respectively)  with index $2$\,.
\qed

The next corollary follows immediately from \eqref{e:fdreq1}.
\begin{Cor}\label{c:n}
If $g$ varies regularly at $\infty$ with index $\alpha \ge 0$. Then, the index $\alpha$ must be in  $[0, 2]$.
\end{Cor}
\proof[Proof of Theorem \ref{t:tail}]
  (i) By Proposition    \ref{prop:Hregvarn} we have that  $\psi$ varies regularly at $0$ with index $2 \wedge \alpha$\,.
    Moreover,  by Proposition \ref{prop:Hregvar2} we have that $H$  varies regularly at $0$  with index $   \alpha/2$\,.
    Noting that  
\begin{align}
\label{e:psiHg}
   \frac{\psi(r^{-1})}{1-\frac{r^{-1}\psi'(r^{-1})}{2\psi(r^{-1})}}=\frac{\psi(r^{-1})^2}{H(r^{-2})}\frac{H(r^{-2})}{g(r^{-1})}, 
\end{align}
    we now use Propositions \ref{prop:Hregvar2} and \ref{prop:tail}  to get
    \begin{align*}
     \lim_{{r\to\infty, \frac{t\psi(r^{-1})}{1-\frac{r^{-1}\psi'(r^{-1})}{2\psi(r^{-1})}}\to 0}}\frac{\P(|X_t|\geq r) }{tg(r^{-1})}&  =
      \lim_{{r\to\infty,  \frac{t\psi(r^{-1})}{1-\frac{r^{-1}\psi'(r^{-1})}{2\psi(r^{-1})}}\to 0}}\frac{\P(|X_t|\geq r) }{tH(r^{-2})}
      \frac{H(r^{-2})}{g(r^{-1})}
      \\
     &=
     \frac{1}{\Gamma(2-\frac{\alpha}{2})}\frac{2^{\alpha-1}\Gamma\left(\frac{d+\alpha}{2}\right)}{\Gamma\left(\frac{d}{2}\right)}.
    \end{align*}
    (ii) The first claim in (ii) is Corollary \ref{c:n}. The proof of the second claim  in (ii)   is analogous to   the proof of (i).
 \qed

\begin{Rem}
If $\psi$ varies regularly at $\infty$, with index $\alpha\in [0,2)$ then 
 it follows from Proposition \ref{prop:hphi-regvar} that  $H$ varies regularly with index $\alpha/2$. 
  Thus Propositions \ref{prop:hphi-regvar} and \ref{prop:tail}  yield
    \begin{align*}
    \lim_{\substack{r\to\infty  (r \to 0, \text{ resp.})\\t\psi(r^{-1})\to 0}}
    \frac{\P(|X_t|\geq r) }{t\psi(r^{-1})}&=\lim_{\substack{r\to\infty  (r \to 0, \text{ resp.})\\ t\frac{\psi(r^{-1})^2}{H(r^{-2})}\to 0}}
    \frac{\P(|X_t|\geq r) }{tH(r^{-2})}\frac{H(r^{-2})}{\psi(r^{-1})}\\
    &=\frac{1}{\Gamma(2-\frac{\alpha}{2})}\frac{(2-\alpha)2^{\alpha-1}\Gamma(\frac{d+\alpha}{2})}{\Gamma(\frac{d}{2})}.    \end{align*}

\end{Rem}

\section{Asymptotic of the transition density  of unimodal L\' evy processes}

Throughout this section we assume that $X$ is a unimodal L\' evy process in $\R^d$, $d\geq 1$ with the characteristic exponent $\xi \to \psi(|\xi|)$, and that 
$p(t,x)$ and $J(x)$ are the transition density and the jumping kernel of $X_t$, respectively.

We will give 
 the precise off-diagonal asymptotic expression of  $p(t,x)$. Recall that $p(t,x)=q(t,|x|)$ where $r \to q(t,r)$ is decreasing.
We also recall that $\phi:(0,\infty)\rightarrow (0,\infty)$ is defined in \eqref{eq:phi}.

\begin{Thm}\label{thm:main2}
  Assume that $H(\lambda)=\phi(\lambda)-\lambda\phi'(\lambda)$ varies regularly at $0$ (at $\infty$, respectively) with index 
 $\gamma\in (0,2)$ and 
$\psi$ varies regularly at $0$ (at $\infty$, respectively)  with index $\gamma_1 \in [0,2]$.
 Then
\begin{align*}
p(t,x)  \sim\frac{\gamma \pi^{-d/2}\Gamma\left(\frac{d}{2}\right)}{\Gamma(2-\gamma)} t|x|^{-d}H(|x|^{-2}) ,
& \quad  |x|\to\infty \,\, (x \to 0, \text{ resp.})
\\
  \text{ and }\quad &
  t\frac{\psi(|x|^{-1})^2}{H(|x|^{-2})}\to 0\,.
\end{align*}
\end{Thm}
\proof
Let $0<a<b$. By Proposition \ref{prop:tail} we have
\begin{align*}
&\lim_{\substack{r\to\infty  (r \to 0, \text{ resp.})\\t\frac{\psi(\sqrt{r}^{-1})^2}{H(r^{-1})}\to 0}}\frac{\P(b^{-1}r\leq |X_t|<a^{-1}r)}{tH(r^{-2})}\\=&\lim_{\substack{r\to\infty  (r \to 0, \text{ resp.})\\t\frac{\psi(\sqrt{r}^{-1})^2}{H(r^{-1})}\to 0}}\left(\frac{\P(|X_t|\geq b^{-1}r)}{tH(r^{-2})}-\frac{\P(|X_t|\geq a^{-1}r)}{tH(r^{-2})}\right)
=\frac{b^{2\gamma}-a^{2\gamma}}{\Gamma(2-\gamma)}.
\end{align*}

Hence, by using 
\[
    \P(|X_t|\geq r)=c_d\int_r^\infty q(t,s)s^{d-1}\,ds,
\]
where $c_d=\frac{2\pi^{d/2}}{\Gamma(d/2)}$,
and the fact that $p(t,x)$ is radially decreasing, for $a<1<b$ we obtain
\begin{align*}
&\limsup_{\substack{|x|\to\infty  (|x| \to 0, \text{ resp.})\\t\frac{\psi(|x|^{-1})^2}{H(|x|^{-2})}\to 0}}\frac{p(t,x)c_dd^{-1}(1-b^{-d})}{t|x|^{-d}H(|x|^{-2})}\\
\leq& \limsup_{\substack{r\to\infty  (r \to 0, \text{ resp.})\\t\frac{\psi(r^{-1})^2}{H(r^{-2})}\to 0}} \frac{\P(|X_t|\geq b^{-1}r)-\P(|X_t|\geq r)}{tH(r^{-2})}\leq \frac{b^{2\gamma}-1}{\Gamma(2-\gamma)}
\end{align*}
and
\begin{align*}
&\liminf_{\substack{|x|\to\infty  (|x| \to 0, \text{ resp.})\\t\frac{\psi(|x|^{-1})^2}{H(|x|^{-2})}\to 0}}\frac{p(t,x)c_dd^{-1}(a^{-d}-1)}{t|x|^{-d}H(|x|^{-2})}\\\geq& \liminf_{\substack{r\to\infty  (r \to 0, \text{ resp.})\\t\frac{\psi(r^{-1})^2}{H(r^{-2})}\to 0}} \frac{\P(|X_t|\geq r)-\P(|X_t|\geq a^{-1}r)}{tH(r^{-2})}\geq \frac{a^{2\gamma}-1}{\Gamma(2-\gamma)}.
\end{align*}
Note that in the case $\gamma=0$ formula holds. 

If $\gamma>0$, we let $a\uparrow 1$ and $b\downarrow 1$ to obtain
\[
    \lim_{\substack{|x|\to\infty  (|x| \to 0, \text{ resp.})\\t\frac{\psi(|x|^{-1})^2}{H(|x|^{-2})}\to 0}}\frac{p(t,x)}{t|x|^{-d}H(|x|^{-2})}=\frac{\gamma \pi^{-d/2}\Gamma\left(\frac{d}{2}\right)}{\Gamma(2-\gamma)}.
\]
\qed

By a consequence of Theorem \ref{thm:main2} we 
prove Theorem \ref{cor:conseq2}
which cover \cite[Theorems 4 and 5]{CGT}  and we further obtain a new result for $\alpha=2$.

\proof[Proof of Theorem \ref{cor:conseq2}]
By Propositions    \ref{prop:Hregvarn} and \ref{prop:Hregvar2} 
 we have that  $\psi$ varies regularly at $0$ (at $\infty$, respectively) with index $2 \wedge \alpha$\,.
  and $H$  varies regularly at $0$ (at $\infty$, respectively) with index $   \alpha/2$.
    By \eqref{e:psiHg}, \eqref{e:P2} now follows immediately from  Theorem \ref{thm:main2} and Proposition \ref{prop:Hregvar2} to get
    \begin{align*}
     \lim_{\substack{|x|\to\infty  (|x| \to 0, \text{ resp.})\\\frac{t\psi(|x|^{-1})}{1-\frac{|x|^{-1}\psi'(|x|^{-1})}{2\psi(|x|^{-1})}}\to 0}}\frac{p(t,x)}{t|x|^{-d}g(|x|^{-1})}
         & =\frac{\alpha\pi^{-d/2}\Gamma(\frac{d}{2})}{\Gamma(2-\frac{\alpha}{2})}\frac{2^{\alpha-1}\Gamma\left(\frac{d+\alpha}{2}\right)}{\Gamma\left(\frac{d}{2}\right)}\nonumber\\
     &=\alpha 2^{\alpha-1}\pi^{-d/2}\frac{\Gamma\left(\frac{d+\alpha}{2}\right)}{\Gamma\left(2-\frac{\alpha}{2}\right)}.
    \end{align*}
\eqref{e:J2} follows from \eqref{e:P2} and the fact that $\lim_{t \downarrow  0} t^{-1} p(t,x)=J(x)$ vaguely on $\R^d\setminus \{0\}$ 
(see \cite[proof of Theorem 6]{CGT} too). 
 \qed

\begin{Rem}

(i)
From the proof of Theorem \ref{cor:conseq2}, we see that  if $g$ varies slowly at $0$ (at $\infty$, respectively) and $\psi$ does not have diffusion part  then
  \begin{align*}
     \lim_{\substack{|x|\to\infty  (|x| \to 0, \text{ resp.})\\\frac{t\psi(|x|^{-1})}{1-\frac{|x|^{-1}\psi'(|x|^{-1})}{2\psi(|x|^{-1})}}\to 0}}\frac{p(t,x)}{t|x|^{-d}g(|x|^{-1})}=0.
    \end{align*}

     (ii) 
    If $\psi$ varies regularly at $0$ (at $\infty$, respectively), with index $\alpha\in [0,2)$ then 
 it follows from Proposition \ref{prop:hphi-regvar} that  $H$ varies regularly at $0$ (at $\infty$, respectively) with index $\alpha/2$.    
    Thus, Theorem \ref{thm:main2} and Proposition \ref{prop:hphi-regvar} yield
    \begin{align*}
    \lim_{\substack{|x|\to\infty  (|x| \to 0, \text{ resp.})\\t\psi(|x|^{-1})\to 0}}
    \frac{p(t,x)}{t|x|^{-d}\psi(|x|^{-1})}&=\lim_{\substack{|x|\to\infty  (|x| \to 0, \text{ resp.})\\ t\frac{\psi(|x|^{-1})^2}{H(|x|^{-2})}\to 0}}
    \frac{p(t,x)}{t|x|^{-d}H(|x|^{-2})}\frac{H(|x|^{-2})}{\psi(|x|^{-1})}\\
    &=\frac{2^{-1}\alpha\pi^{-d/2}\Gamma(\frac{d}{2})}{\Gamma(2-\frac{\alpha}{2})}\frac{(2-\alpha)2^{\alpha-1}\Gamma(\frac{d+\alpha}{2})}{\Gamma(\frac{d}{2})}\\
    &=\alpha 2^{\alpha-1}\pi^{-d/2}\frac{\Gamma(\frac{d+\alpha}{2})}{\Gamma(1-\frac{\alpha}{2})}.
    \end{align*}
\end{Rem}

We apply Theorem \ref{cor:conseq2} to subordinate Brownian motion $X=(X_t)_{t\geq 0}$. It is an isotropic L\' evy process in  $\R^d$ $(d\geq 1)$ defined by
\[
    X_t=B_{T_t},\quad t\geq 0\,,
\]
where $B=(B_t)_{t\geq 0}$ is the $d$-dimensional Brownian motion with the transition density \[
p^{(2)}(t,x)=(4\pi t)^{-d/2}e^{-\frac{|x|^2}{4t}},\quad t>0,\,\,x\in \R^d\,.
\]and $T=(T_t)_{t\geq 0}$ is an independent subordinator whose the Laplace exponent is $\varphi$. 

It is known that the characteristic exponent of $X$ is given by $\psi(\xi)=\varphi(|\xi|^2)$ and $X$ has a transition transition density given by
\begin{equation}\label{eq:sbm-td}
p(t,x)=\int_{(0,\infty)}(4\pi s)^{-d/2}e^{-\frac{|x|^2}{4s}}\P(T_t\in ds)\,.
\end{equation}
Thus, subordinate Brownian motion is a subclass of unimodal L\' evy processes and the next corollary follows directly from 
Theorem \ref{cor:conseq2}.

\begin{Cor}\label{cor:conseq3}
$X_t=B_{T_t}$ is a subordinate Brownian motion and $\varphi$ is the Laplace exponent of the subordinator $T_t$.
%(i) If $\varphi$ varies regularly at $0$ (at $\infty$, respectively) with index $\gamma\in (0, 1)$, then 
%\begin{align*}
%p(t,x)  &\sim \gamma 4^{\gamma}\pi^{-d/2}\frac{\Gamma(\frac{d}{2} +\gamma )}{\Gamma(1-\gamma)}  t|x|^{-d}\varphi(|x|^{-2})
%\\&=\gamma 4^{\gamma}\pi^{-d/2-1}\sin\left(\gamma \pi\right)\Gamma\left(\gamma\right)\Gamma\left(\frac{d}{2} +\gamma\right)  t|x|^{-d}\varphi(|x|^{-2}),
%\\
% &\quad  |x|\to\infty \,\, (x \to 0, \text{ resp.})
% \quad \text{ and }  \quad
%  t \varphi(|x|^{-2})\to 0\,.
%\end{align*}
%
%\noindent
%(ii) 
Suppose  $H(\lambda)=\varphi(\lambda)-\lambda\varphi'(\lambda)$ varies regularly at $0$ (at $\infty$, respectively) with index $\gamma\in (0, 2)$. 
We further assume that $\varphi$ does not have drift part if $H$ varies regularly at $\infty$ with index $\gamma\in (0, 1)$.
 Then 
\begin{align}
p(t,x)  &\sim
\gamma 4^{\gamma}\pi^{-d/2}\frac{\Gamma\left(\frac{d}{2} +\gamma\right)}{\Gamma\left(2-\gamma\right)}t|x|^{-d}H(|x|^{-2})
\nonumber\\
 &\quad  |x|\to\infty \,\, (x \to 0, \text{ resp.})
 \quad \text{ and }  \quad
 \frac{t\varphi(|x|^{-2})^2}{H(|x|^{-2})}\to 0\,.    \label{c:enw}
\end{align}
\end{Cor}

\begin{Ex}
Let us consider a Bernstein function $\varphi(\lambda)=\lambda\log(1+\frac{1}{\lambda})$. 
(see table entry No. 27 on page 316 in \cite{SSV}). 
Then $\varphi$ varies regularly at $0$ with index $1$ and 
we have $H(\lambda)=\lambda/(\lambda+1)$ so that 
\[
    \lim_{\lambda \downarrow 0}\frac{H(\lambda)}{\lambda}=1\quad \text{and }\quad  \frac{\varphi(\lambda)}{H(\lambda)}\sim \log(1+\frac{1}{\lambda}),\,\,\,\,\lambda\downarrow 0.
\]
By \eqref{c:enw} we have the asymptotic formula as 
\[
     \lim_{\substack{|x|\to\infty\\t|x|^{-2}\log(1+|x|)^2\to 0}}\frac{p(t,x)}{t|x|^{-d-2}}=4\pi^{-d/2}\Gamma\left(\frac{d+2}{2}\right).
\]
\end{Ex}

\noindent
{\bf Acknowledgements:} The first named author is grateful to 
%Juhak Bae
Joohak Bae for
reading the earlier draft of this paper and giving helpful comments.
We thank the referees for helpful comments on the first version of this paper.

\bibliography{myrefs}
\end{document}